\documentclass[11pt,reqno,toc=bibliography,]{scrartcl}
\usepackage[margin=2cm
]{geometry}
\makeatletter
\DeclareOldFontCommand{\rm}{\normalfont\rmfamily}{\mathrm}
\DeclareOldFontCommand{\sf}{\normalfont\sffamily}{\mathsf}
\DeclareOldFontCommand{\tt}{\normalfont\ttfamily}{\mathtt}
\DeclareOldFontCommand{\bf}{\normalfont\bfseries}{\mathbf}
\DeclareOldFontCommand{\it}{\normalfont\itshape}{\mathit}
\DeclareOldFontCommand{\sl}{\normalfont\slshape}{\@nomath\sl}
\DeclareOldFontCommand{\sc}{\normalfont\scshape}{\@nomath\sc}
\makeatother

\usepackage[utf8]{inputenc}

\usepackage[english]{babel}

\usepackage{physics}
\usepackage{braket}
\usepackage{dsfont}
\usepackage[table,xcdraw]{xcolor}
\usepackage{multirow}
\usepackage{graphicx}
\usepackage{tikz}
\usepackage{amsmath}
\usepackage{amssymb}
\usepackage{amsthm}
\usepackage{mathtools}
\usetikzlibrary{matrix}
\usepackage{appendix}
\usepackage{lipsum}
\usepackage{hyperref}
\usepackage{mathrsfs}

\usepackage{xcolor}
\usepackage{bbm} 

\newtheoremstyle{break}
  {5pt}{5pt}%
  {\itshape}{}%
  {\bfseries}{}%
  {\newline}{}%
\theoremstyle{break}

\usepackage[color]{changebar}

\usepackage{tocbasic}

\newtheorem{definition}{Definition}[section]
\newtheorem{theorem}[definition]{Theorem}

\newtheorem{lemma}[definition]{Lemma}
\newtheorem{corollary}[definition]{Corollary}
\newtheorem{remark}[definition]{Remark}

\newcommand{\refdef}[1]{\textbf{Definition~\ref{#1}}}
\newcommand{\reflem}[1]{\textbf{Lemma~\ref{#1}}}

\newcommand{\refcor}[1]{\textbf{Corollary~\ref{#1}}}

\newcommand{\refthm}[1]{\textbf{Theorem~\ref{#1}}}

\renewcommand{\d}{\,\mathrm{d}}
\renewcommand{\i}{\mathrm{i}}
\newcommand{\e}{\mathrm{e}}

\renewcommand{\d}{\,\mathrm{d}}
\renewcommand{\dd}{\mathrm{d}}

\renewcommand{\Re}{\mathrm{Re}}

\newcommand{\note}[1]{} 
\newcommand{\N}{\mathbbm{N}}

\newcommand{\R}{\mathbbm{R}}

\newcommand{\C}{\mathbbm{C}}

\newcommand{\B}{\mathbbm{B}}
\renewcommand{\H}{\mathcal{H}}
\newcommand{\foral}{\forall\,}

\newcommand{\D}{\mathrm{D}}
\renewcommand{\epsilon}{\varepsilon}

\usepackage{tcolorbox}

\newtcolorbox[auto counter, number within=section]{mybox}[2][]{
  colframe=black,  
  colback=white,   
  fonttitle=\bfseries, 
  coltitle=black,  
  colbacktitle=white, 
  enhanced,
  sharp corners, 
  boxed title style={colframe=white}, 
  title={#2}, 
}

\usepackage[normalem]{ulem} 
\usepackage{xcolor}

\setkomafont{author}{\sffamily}
\setkomafont{title}{\sffamily}
\makeatletter
\patchcmd{\@maketitle}{\huge}{\LARGE}{}{}
\makeatother

\RedeclareSectionCommand[tocbeforeskip=1ex plus 1pt minus 1pt]{section}
\makeatletter
\renewcommand*{\@pnumwidth}{1em}
\makeatother

\deftocheading{toc}{} 

\title{\vspace*{-22pt} On local well-posedness for the nonlinear Schrödinger equation with general power nonlinearity}
\subtitle{Two Different Approaches}

\usepackage{authblk}

\author{Lucia Arens,\ Marius Gritl}

\affil{\small Department of Mathematics, Technische Universit\"at M\"unchen, Germany\\
{\footnotesize\texttt {\{lucia.arens, m.gritl\}@tum.de}}}

\date{\vspace*{-30pt}\small \vspace*{-12pt} }

\begin{document}  

\maketitle

\tableofcontents

\begin{abstract}\vspace*{5pt} \noindent{\bf Abstract: }\small
The nonlinear Schrödinger equation plays a fundamental role in mathematical physics, particularly in the study of quantum mechanics and Bose-Einstein condensation. This paper explores two distinct approaches to establishing the local well-posedness of solutions: the semigroup theory ansatz and the Strichartz estimates ansatz. Semigroup theory provides a general and elegant framework rooted in functional analysis, allowing for the interpretation of the time evolution of solutions as operator semigroups. Strichartz estimates, developed specifically for dispersive equations, offer an alternative technique based on refined space-time estimates and fixed-point arguments. We systematically analyze and compare both approaches and apply them to nonlinear Schrödinger equations where the nonlinearity is given by $F(u)=\lambda|u|^p u$ for some $\lambda \in \R$. So our results extend beyond the physically relevant case $p=2$. 
\end{abstract}

\section{Introduction}

We consider the Cauchy problem for the nonlinear Schrödinger equation (NSE):
\begin{equation}
    \i \partial_t u + \Delta u = \lambda |u|^p u, \quad u(0, x) = u^0(x),
    \label{eq:NLS}
\end{equation}
with $u: \R \times \R^d \to \C$, nonlinearity parameter $\lambda \in \R$, $p > 0$ and nonlinearity $F(u) = \lambda |u|^pu$. The initial datum $u^0$ is assumed to belong to the space $L^2(\R^d)$.\\

For $\lambda = 0$ the NSE reduces to the well-known Schrödinger equation that characterizes the time evolution of a free particle in the setup of a quantum system. The nonlinear Schrödinger equation for $\lambda \in \R \backslash{ \{ 0 \}}$ arises in various physical settings including fluid dynamics or nonlinear optics, as well as in condensed matter physics, for instance in describing the Bose-Einstein condensation phenomenon or even in mathematical finance. In the context of Bose-Einstein condensation, the NSE is known as the Gross-Pitaevskii equation for $p=2$. The nonlinearity corresponds to interaction effects in the condensate \cite{BEC}.\\
For water waves, the nonlinear Schrödinger equation describes the envelope of two-dimensional gravity waves in water of finite depth. In this context, the nonlinearity parameter $\lambda$ depends on the relative water depth. If the water depth is large compared to the wave length of the water waves, $\lambda$ is negative and envelope solitons may occur. For shallow water $\lambda$ is positive and wave groups with envelope solitons do not exist, see \cite{water_waves} for more details.\\
In optical fiber engineering, the nonlinear Schrödinger equation models wave propagation through a nonlinear medium. The nonlinearity parameter $\lambda$ corresponds to the magnitude of the nonlinear Kerr effect. The defocusing NSE with $\lambda > 0$ has solutions that tend to spread out and remain well-behaved globally since dispersive effects dominate. The focusing NSE with $\lambda < 0$ exhibits finite-time blowup phenomena. It allows for solutions that are localized in space and have spatial attenuation towards infinity or nonlinear waves in which energy concentrates in a localized and oscillatory fashion. The study of blow-up structures remains an ongoing challenge \cite{Blowup}.\\

We will use semigroup theory and Strichartz estimates to analyze the local behavior of solutions. Firstly, the time evolution of solutions is represented by an operator semigroup. This allows us to use the extensive results for semigroups from operator theory, first and foremost \textsc{Stone's theorem}. \textsc{Picard-Lindelöf's Theorem} famously provides a set of conditions under which an initial value problem of an ordinary differential equation has a unique local and global solution. A similar strategy will also prove successful in the context of the NSE. Strichartz estimates provide the necessary norm estimates of solutions in mixed Lebesgue spaces, in particular of the free Schrödinger propagator, and ensure that the nonlinear growth is sufficiently controlled. Strichartz estimates were first introduced in 1977 to study the decay of solutions of the wave equation \cite{Strichartzwave}. A local existence result in a suitable Strichartz space $X$ is derived using \textsc{Banach's fixed-point theorem}. 

\section{Preliminaries}
This section contains some basic definitions and theorems that will be used in the upcoming paper. Many statements are given without proof.

\subsection{Fourier transform}
\label{chap:propagator}

\begin{definition}[Fourier transform]
    If $u \in L^1(\R^d)$, we define its \textbf{Fourier transform}
    \begin{align*}
        \hat{u}(k)=(2\pi)^{-\frac{d}{2}} \int_{\R^d} \e^{-\i k\cdot x}u(x) \d x, \quad k \in \R^d.
    \end{align*}
 Its \textbf{inverse Fourier transform} is defined by
    \begin{align*}
        \check{u}(k)=(2\pi)^{-\frac{d}{2}} \int_{\R^d} \e^{\i k\cdot x}u(x) \d x, \quad k \in \R^d.
    \end{align*}
\end{definition} \label{Fourier}
Since $L^1(\R^d) \cap L^2(\R^d)$ is dense in $L^2(\R^d)$ there is a unique continuous extension of the Fourier transform to $L^2(\R^d)$.

\begin{lemma}[Hausdorff-Young Inequality]
    Let $p \in [1,2]$ and $(p,p')$ be conjugated exponents. Then the following holds for all $f \in L^p(\R^d)$:
    \begin{align}
        \Vert \hat{f} \Vert_{L^{p'}(\R^d)} \leq \Vert f \Vert_{L^p(\R^d)}.\label{hyineq}
    \end{align}
\end{lemma}

In various situations in quantum mechanics one has to define a function of an operator, i.e. an operator-valued function. One can do this via various types of functional calculus; we will do this later. But we want to present procedure for the Laplacian at this point by using the Fourier transform: we define the operator exponential of the Laplacian which occurs in the Schrödinger equation.\\

We consider the linear Schrödinger equation 
\begin{align}
 \begin{cases}
    \begin{aligned}
        \i \partial_t u(t,x) + \Delta u(t,x) & = 0,  \quad &\text{on } \R \times \R^d, \\
        u(0,x) & = u^0(x),  \quad &\text{on } \R^d.
    \end{aligned}
\end{cases}\label{eq:linSE}
\end{align}
One can formally denote the solution as $u(t,x)=\exp(\i t \Delta) u^0(x)$, where $\exp(\i t \Delta)$ is called the \textit{free Schrödinger propagator}. We can define this object by using the Fourier transform in a precise way.

\begin{lemma} 
Let $f\in L^2(\R^d)$, then \begin{align*}
        \mathscr{F}(\Delta f)(k)= - |k|^2 \hat{f}(k), 
    \end{align*} 
    where $\mathscr{F}$ denotes the Fourier transform on $L^2(\R^d)$.
\end{lemma}

Using this elementary fact, one can compute the Fourier transform of the entire Schrödinger equation to obtain an ODE in the momentum space:
\begin{align*}
 \begin{cases}
    \begin{aligned}
        \i \partial_t \hat{u}(t,k)-|k|^2 \hat{u}(t,k)&=0,  \quad &\text{on } \R \times \R^d, \\
        \hat{u}(0,k) & = \widehat{u^0}(k),  \quad &\text{on } \R^d.
    \end{aligned}
\end{cases}
\end{align*}
Clearly, the solution is given by \begin{align*}
    \hat{u}(t,k)=\e^{-\i t |k|^2} \widehat{u^0}(k).
\end{align*}
Finally, we invoke the inverse Fourier transform to obtain the solution $u$ in the original space. This motivates the following definition.

\begin{definition}
    The \textbf{free Schrödinger propagator} $t \mapsto \exp(\i t \Delta)$ is implicitly defined via \begin{align*}
        \exp(\i t \Delta) u^0(x) = \mathscr{F}^{-1}\big( \e^{-\i t |k|^2} \widehat{u^0}(k)\big), \quad t \in \R,
    \end{align*} where $u^0\in L^2(\R^d)$ denotes the initial condition in \eqref{eq:linSE}.\label{propagator}
\end{definition} 

\subsection{Sobolev spaces}

\begin{definition}
Let $\Omega \subseteq \R^d$ be open and non-empty. Let $\alpha \in \N_0^d$ be a multi-index. A function $u\in L_\mathrm{loc}^1(\Omega)$ has the \textbf{weak derivative} $\D^{\alpha}u \in L_\mathrm{loc}^1(\Omega)$ if
\begin{align*}
\int_{\Omega}(\D^{\alpha} u) \phi \d x = (-1)^{|\alpha|} \int_{\Omega} u(\D^{\alpha} \phi) \d x
\end{align*}
holds for all $\phi \in \mathscr{C}_\mathrm{c}^{\infty}(\Omega)$.
\end{definition}

\begin{definition}
Let $\Omega \subseteq \R^d$ be open and non-empty. Let $p \in [1,\infty]$ and $k \in \N$. Then we define the \textbf{Sobolev space} $W^{k,p}(\Omega)$ by
\begin{align*}
    W^{k,p}(\Omega)=\{ u \in L^p(\Omega) : \D^{\alpha}u \in L^p(\Omega) \text{ exists for all } \alpha \in \N_0^d, |\alpha| \leq k \}.
\end{align*}
$H^k(\Omega) := W^{k,2}(\Omega)$ is a Hilbert space with
\begin{align*}
    \langle u,v \rangle_{H^k(\Omega)} := \sum_{|\alpha| \leq k} \langle \D^{\alpha}u, \D^{\alpha}v \rangle_{L^2(\Omega)}.
\end{align*}
 Alternatively, the space  $H^k(\Omega)$ can be defined using the Fourier transform, namely
\begin{align*}
    H^k(\Omega) = \left\{ u \in L^2(\Omega) : \left(1+ | \eta |^2 \right)^\frac{k}{2} \hat{u}(\eta) \in L^2(\R^d) \right\}.
\end{align*}
One can define the norm
\begin{align*}
    \Vert u \Vert_{H^k(\Omega)} = \left( \int_{\Omega} \left(1+ | \eta |^2 \right)^k \vert \hat{u}(\eta) \vert^2 \d \eta \right)^\frac{1}{2},
\end{align*}
which is equivalent to the norm induced by the above defined scalar product.
\end{definition}

The Sobolev spaces $W^{k,p}(\Omega)$ combine the concepts of weak differentiability and Lebesgue norms. Some proofs will use sobolev embeddings.
The following lemma can be found in \cite[Corollary 9.13]{Brezis}.

\begin{lemma}[Sobolev embedding theorem]
Let $m \in \N$ and $p \in [1,\infty)$. We have the following continuous embeddings:
\begin{enumerate}
    \item[(i)] $W^{m,p}(\R^d) \subseteq L^q  (\R^d)$, where $\frac1q=\frac1p-\frac{m}{d}$ if $\frac1p - \frac{m}{N} >0$,
    \item[(ii)]\label{lemmii} $W^{m,p}(\R^d) \subseteq L^q  (\R^d)$ for all $q \in [p,\infty)$ if $\frac1p-\frac{m}{d}=0$,
    \item[(iii)] $W^{m,p}(\R^d) \subseteq L^\infty  (\R^d)$ if $\frac1p-\frac{m}{N}<0$.
\end{enumerate}\label{lemsobemb}
\end{lemma}

\begin{remark}
    At this point, we would like to note two specific embeddings that will be required later, namely
    \begin{align}
        H^2(\R^2) \subseteq L^\infty(\R^2) \label{sobemb1}
      \end{align}  
     and   
    \begin{align}
        H^1(\R^2) \subseteq L^r(\R^2) \, \text{ for } 2 \leq r  <\infty. \label{sobemb2}
    \end{align}
    The first one follows by setting $m=2$, $p=2$ and $d=2$. We compute $\frac1q=\frac12 -\frac{2}{2}=-\frac12 <0$, so the result follows by (iii) in \reflem{lemsobemb}. The second embedding follows by taking $m=1$, $p=2$ and $d=2$. Since $\frac1q = \frac12-\frac12=0$ we can apply (ii).
\end{remark}

\subsection{Self-adjointness of operators on Hilbert spaces}

The following concepts are ubiquitous in operator theory, for more details refer to \cite{Simon}.

\begin{definition}[Notions for operators on Hilbert spaces]
Let $\H$, $\Tilde{\H}$ be Hilbert spaces and $T:\mathscr{D}_T\to \Tilde{\H}$ be a linear operator where $\mathscr{D}_T \subseteq \H$ is a subspace of $\H$, called the \textbf{domain} of $T$.
\begin{itemize}
    \item $T$ is called \textbf{densely defined} if $\overline{\mathscr{D}_T}=\H$.
    \item The \textbf{adjoint} $T^*$ of a densely defined operator $T$ is defined by the condition \begin{align*}
        \foral \psi \in \mathscr{D}_T\, \foral \varphi \in \mathscr{D}_{T^*} : \langle T^* \varphi, \psi\rangle = \langle \varphi, T \psi\rangle,
    \end{align*} 
    where $\mathscr{D}_{T^*}:=\{\varphi \in \Tilde{\H}| \mathscr{D}_T \ni \psi \mapsto \langle \varphi, T \psi\rangle \text{ is bounded}\}$.
    \item A densely defined operator $T$ is called \textbf{symmetric} if $\langle T \varphi, \psi \rangle = \langle \varphi, T \psi \rangle$ holds for all $\varphi, \psi \in \mathscr{D}_T$.
    \item $T$ is called \textbf{self-adjoint} if $T=T^*$, that is, iff $T$ is symmetric and $\mathscr{D}_T=\mathscr{D}_{T^*}$ holds true.
    \item $T$ is called \textbf{closed} if the graph $\Gamma(T):=\{\langle \varphi, T \varphi\rangle | \varphi \in \mathscr{D}_T\}$ is closed w.r.t. the product topology on $\H \times \Tilde\H$.
    \item An operator $S$ on $\H$ is called an \textbf{extension} of $T$ and we write $S \supseteq T$, if $\mathscr{D}_S \supseteq \mathscr{D}_T$ and $S\varphi = T \varphi$ holds for all $\varphi \in \mathscr{D}_T$. 
    \item $T$ is \textbf{closable} if it has an closed extension. The smallest closed extension, which is unique, is called the \textbf{closure} $\overline{T}$.
    \item $T$ is called \textbf{essentially self-adjoint}, if its closure $\overline{T}$ is self-adjoint.
\end{itemize}
\end{definition}

\begin{lemma}
    Let $T: \mathscr{D}_T \to \H$ be a densely defined operator on $\H$. Then $T^*$ is closed.\label{lemmaadjointclosed}
\end{lemma}

\begin{proof}
    Let $(\varphi_n, \eta_n) \in \Gamma(T^*)$ be a convergent sequence with limit $(\varphi, \eta) \in \H \times \H$. Then for all $\psi \in \mathscr{D}_T$, the expressions on the left and right side in \begin{align*}
        \langle \varphi_n, T \psi \rangle = \langle \eta_n, \psi \rangle
    \end{align*}
    converge to \begin{align*}
        \langle \varphi, T \psi \rangle = \langle \eta, \psi \rangle ,
    \end{align*}
    respectively. From that we can conclude $\varphi \in \mathscr{D}_{T^*}$ and $T^* \varphi = \eta$, so $T^*$ is closed.
\end{proof}

\begin{lemma}
    Let $T: \mathscr{D}_T \to \H$ be a symmetric operator on $\H$. Then $T \subseteq T^*$.
\end{lemma}

\begin{proof}
    Let $\varphi \in \mathscr{D}_T$. The map $\psi \mapsto \langle \varphi, T \psi \rangle$ is bounded, since the following holds for every $\psi \in \mathscr{D}_T$: \begin{align*}
        |\langle \varphi, T \psi\rangle|=|\langle T\varphi, \psi\rangle| \leq \|T \varphi\|\|\psi\|.
    \end{align*}
    Hence, $\varphi$ belongs to $\mathscr{D}_{T^*}$. By assumption, we have $T^* \varphi = T \varphi$ for all $\varphi \in \mathscr{D}_T$.
\end{proof}

\begin{theorem}[Basic criterion for self-adjointness]
Let $T$ be a symmetric operator on a Hilbert space $\H$. Then the following statements are equivalent:
\begin{enumerate}
    \item[(i)] $T$ is self-adjoint,
    \item[(ii)] $T$ is closed and $\ker(T^* \pm \i \mathds1)=\{0\}$,
    \item[(iii)] $\mathrm{im}(T \pm \i \mathds1)=\H$.
\end{enumerate}\label{criterionsa}
\end{theorem}

\begin{proof}
\begin{enumerate}
    \item[(i)] $\Longrightarrow$ (ii): Thanks to \reflem{lemmaadjointclosed}, $T=T^*$ is closed. \\Let w.l.o.g. $\varphi \in \ker(T^* + \i \mathds{1})\backslash\{0\}$. Then one has \begin{align*}
        \i \langle \varphi, \varphi \rangle = \langle - \i \varphi, \varphi \rangle = \langle T^* \varphi, \varphi \rangle = \langle \varphi, T^*\varphi\rangle = \langle \varphi, -\i\varphi\rangle=-\i\langle \varphi, \varphi\rangle.
    \end{align*}
    Hence, $\varphi \equiv 0$.
    \item[(ii)] $\Longrightarrow$ (iii): From \begin{align*}
        \varphi \in \ker(T^*+\i \mathds{1}) \Longleftrightarrow \foral \psi \in \mathscr{D}_T: \langle (T^*+ \i \mathds1)\varphi, \psi \rangle =0 \Longleftrightarrow \foral \psi \in \mathscr{D}_T: \langle \varphi, (T-\i \mathds1)\psi\rangle =0,
    \end{align*}
    we can conclude $\ker(T^* + \i \mathds1)=\mathrm{im}(T-\i \mathds1)^\perp$ and therefore, by assumption, \begin{align*}
        \H = \ker(T^* + \i \mathds1)^\perp = \left(\mathrm{im}(T-\i \mathds1)^\perp\right)^\perp = \overline{\mathrm{im}(T-\i \mathds1)}.
    \end{align*}
    Exploiting this denseness, we choose an $\eta \in \H$ s.t. $(T-\i \mathds1)\varphi_n \to \eta$ for $(\varphi_n)_{n\in\N}\subseteq \H$.\\ The following inequality holds true for all $\psi \in \mathscr{D}_T$: \begin{align*}
        \|(T-\i \mathds1)\psi\|^2 = \langle (T-\i \mathds1)\psi, (T-\i \mathds1)\psi\rangle = \|T\psi\|^2 + \|\psi\|^2 + \i (\langle \psi, T \psi\rangle - \langle T \psi, \psi\rangle) \geq \|\psi\|^2.
    \end{align*}
    This allows us to transfer the Cauchy property from $((T-\i \mathds1)\varphi_n)_{n\in\N}$ to $(\varphi_n)_{n\in\N}$. Hence, there is a $\varphi \in \H$ s.t. $\varphi_n \to \varphi$.\\
    Since $T$ is closed, the limit $\varphi$ is an element of $\mathscr{D}_T$ and $T \varphi=\eta + \i \varphi$. This shows $\eta \in \mathrm{im}(T-\i \mathds1)$; so the surjectivity $\mathrm{im}(T-\i \mathds1)=\H$ is proven.
    \item[(iii)] $\Longrightarrow$ (i): Let $\varphi \in \mathscr{D}_{T^*}$. By assumption, there is a $\psi \in \mathscr{D}_T$ with \begin{align*}
        (T^*-\i \mathds1)\varphi = (T-\i \mathds1)\psi.
    \end{align*}
    Since $T \subseteq T^*$, we thus have \begin{align*}
        (T^*-\i \mathds1)(\varphi-\psi)=0. 
    \end{align*}
    Then, for all $\eta \in \mathscr{D}_T$:
    \begin{align*}
        0 = \langle \eta, (T^*-\i \mathds1)(\varphi-\psi)\rangle=\langle(T+\i \mathds1)\eta, \varphi-\psi\rangle,
    \end{align*}
    so $\varphi = \psi \in \mathscr{D}_T$ because of $\mathrm{im}(T+\i \mathds1)=\H$.
\end{enumerate}
\end{proof}

\begin{corollary}[Basic criterion for essential self-adjointness]
    Let $T$ be a symmetric operator on a Hilbert space $\H$. Then the following statements are equivalent:
\begin{enumerate}
    \item[(i)] $T$ is essentially self-adjoint,
    \item[(ii)] $\ker(T^* \pm \i \mathds1)=\{0\}$,
    \item[(iii)] $\overline{\mathrm{im}(T \pm \i \mathds1)}=\H$.
\end{enumerate}\label{criterionesa}
\end{corollary}


\subsection{Continuous one-parameter groups and Stone's theorem}

Continuous one-parameter unitary groups play an important role in quantum mechanics as they are useful to describe the time evolution of quantum systems. This chapter is based on \cite[Chapter VIII.4]{Simon}.

\begin{definition}
    A \textbf{strongly continuous one-parameter unitary group} $\{U(t)\}_{t\in \R}$ satisfies
    \begin{enumerate}
        \item \textbf{unitarity}: $U(t)^*=U(t)^{-1}$ for all $t \in \R$,
        \item \textbf{group property}: $U(t)U(s) = U(s)U(t)=U(t+s)$ for all $t,s \in \R$ and $U(0) =\mathds1$,
        \item \textbf{strong continuity}: $\lim_{t \to t_0} U(t)\varphi=U(t_0)\varphi$ for all $t_0 \in \R$ and $\varphi \in \H$.
    \end{enumerate}
\end{definition}

\begin{lemma}
    Let $A: \mathscr{D}_A \to \H$ be a self-adjoint operator and define $U(t)=\exp(\i tA)$. Then $\{U(t)\}_{t\in\R}$ is a strongly continuous one-parameter unitary group. Furthermore, one has \begin{enumerate}
        \item For $\psi \in \mathscr{D}_A$, \begin{align}
            \lim_{t \to 0} \frac{U(t)\psi - \psi}{t} = \i A \psi, \label{lemmastone1}
        \end{align}
        \item If $\lim_{t \to 0} \frac{U(t)\psi - \psi}{t}$ exists, then \begin{align}\psi \in \mathscr{D}_A.\label{lemmastone2}\end{align}
\end{enumerate}\label{lemmastone}
\end{lemma}

\begin{remark}
The operator $A$ is called the \textbf{infinitesimal generator} of $\{U(t)\}_{t \in \R}$.
\end{remark}

\begin{proof}
    Unitarity and the group property are immediate consequences of the functional calculus. Let $P$ be the unique spectral measure associated with the self-adjoint operator $A$. Then the spectral theorem states that \begin{align*}
        f(A)= \int_{\R} f(\lambda) P(\dd \lambda),
    \end{align*}
    where $f:\R \to \C$ is a measurable function. 
    \\ 
    Inserting $f:\R \to \C, f(x)=\e^{\i t x}$, we can conclude
    \begin{align*}
        U(t)^* = \left( \int_{\R} \e^{\i t \lambda} P(\dd \lambda)\right)^* = \int_{\R} \e^{-\i t \lambda} P(\dd \lambda) = U(-t),
    \end{align*}
    and \begin{align*}
        U(t)U(s)= \exp(\i t A) \exp(\i s A) = \int_{\R} \e^{\i t \lambda} \e^{\i s \lambda} P(\dd\lambda) = \int_{\R} \e^{\i(t+s)\lambda} P(\dd\lambda) = U(t+s).
    \end{align*}
    In order to prove strong continuity, we use \begin{align*}
        \|\exp(\i tA) \varphi - \varphi\|^2 = \int_{\R} |\e^{\i t \lambda}-1|^2 \d \langle P \varphi, \varphi\rangle.
    \end{align*}
    Since $t\mapsto |\e^{\i t \lambda}-1|^2$ is bounded by $2$ and satisfies \begin{align*}
        \lim_{t \to 0}|\e^{\i t \lambda} -1|^2 =0
    \end{align*}
    for all $\lambda \in \R$, we end up with \begin{align*}
        \lim_{t \to 0} \|U(t) \varphi - \varphi\|=0.
    \end{align*}
    This ensures strong continuity at $t=0$. Thanks to the group property, we can conclude strong continuity of $t \mapsto U(t)$ everywhere. \\
    \eqref{lemmastone1} can be proven by the same arguments. One exploits the \textsc{dominated convergence theorem} again and uses $|\e^{\i x}-1|\leq |x|$.\\ 
    To prove \eqref{lemmastone2}, we use \eqref{lemmastone1}: the infinitesimal generator $A$ can be obtained by \begin{align*}
        A \psi = - \i \lim_{t \to 0} \frac{U(t) \psi - \psi}{t}.
    \end{align*}
    Since we assume that this limit exists, it follows that $\psi \in \mathscr{D}_A$.
\end{proof}

Stone's theorem establishes a one-to-one correspondence between self-adjoint operators and strongly continuous one-parameter unitary groups.

\begin{theorem}[\textsc{Stone's theorem}]
Let $\{U(t)\}_{t \in \R}$ be a strongly continuous one-parameter unitary group. Then there exists a unique (possibly unbounded) self-adjoint operator $A: \mathscr{D}^\mathrm{Stone}_A \rightarrow \H$ such that
\begin{align*}
    U(t)=\exp(\i tA)
\end{align*}
holds for all $t \in \R$. The domain of $A$ is defined by
\begin{align*}
    \mathscr{D}^\mathrm{Stone}_A := \left\{ \phi  \in \H : \lim_{\epsilon \to 0} \frac{-\i(U(\epsilon)\phi - \phi)}{\epsilon}\, exists \right\}.
\end{align*}
\end{theorem}

\begin{proof}

\begin{enumerate}
    \item \textit{Define a proper subspace $\mathscr{D}$ that is dense in $\H$:}
    We define \begin{align*}
        \varphi_f := \int_\R f(t) U(t) \varphi \d t
    \end{align*}
    for all $\varphi \in \H$ and some $f \in \mathscr{C}_\text{c}^\infty(\R)$.
    Let $\mathscr{D}$ be the set of linear combinations of all such $\varphi_f$. We define the approximate identity to be the sequence $\{j_\epsilon\}$ of functions with support in $[-\epsilon,\epsilon] \subseteq \R$ for $\epsilon>0$, where \begin{align*}
        j_\epsilon(x) \geq 0  \text{ for all } x \in \R \quad \text{and} \quad \int_\R j_\epsilon(x) \d x =1.
    \end{align*} Then \begin{align*}
        \|\varphi_{j_\epsilon}-\varphi\| &= \left\| \int_\R j_\epsilon(t) (U(t)\varphi-\varphi) \d t \right\|\\
        &\leq \int_\R j_\epsilon(t) \d t \sup_{t \in[-\epsilon,\epsilon]} \|U(t)\varphi-\varphi\|.
    \end{align*} 
    This shows -- thanks to the strong continuity of $U(t)$ -- that $\mathscr{D}$ is dense in $\H$.
    \item \textit{Differentiate $U(t)$ at $t=0$ to obtain $A$:}
    By identity \eqref{lemmastone1} one can suggest that we are be able to obtain $A$ by differentiating $U(t)$ at $t=0$. For $\varphi_f \in \mathscr{D}$, we compute \begin{align*}
        \frac{U(s)-\mathds1}{s} \varphi_f &= \frac{U(s)-\mathds1}{s} \int_\R f(t) U(t) \varphi \d t \\
        &= \int_\R f(t) \frac{U(s+t)-U(t)}{s} \varphi \d t \\ 
        &=\underbrace{\int_\R f(t) \frac{U(s+t)}{s} \varphi \d t}_{\text{substitute } \tau=s+t} - \underbrace{\int_\R f(t) \frac{U(t)}{s} \varphi \d t}_{\text{substitute } \tau=t} \\
        &= \int_\R f(\tau-s) \frac{U(\tau)}{s} \varphi \d \tau - \int_\R f(\tau) \frac{U(\tau)}{s} \varphi \d \tau \displaybreak[1] \\
        &= \int_\R \underbrace{\frac{f(\tau-s)-f(\tau)}{s}}_{\stackrel{s \to 0}{\rightarrow} -f'(\tau) \text{ uniformly}} U(\tau) \varphi \d \tau\\
        &\hspace{-0.1cm}\stackrel{s \to 0}{\rightarrow} - \int_\R f'(\tau) U(\tau) \varphi \d t \\
        &= \varphi_{-f'}.
        \end{align*} 
Since for all $\varphi_f \in \mathscr{D}$ we have shown that the limit exists, this implies that $\varphi_f \in \mathscr{D}_A^{\mathrm{Stone}}$ and therefore $\mathscr{D} \subseteq \mathscr{D}_A^{\mathrm{Stone}}$. 
Inspired by \eqref{lemmastone1}, we define $A \varphi_f = -\i\varphi_{-f'}$. This operator is a restriction of the operator defined by the derivative at $t = 0 $ as stated in the theorem.

    \item \textit{$A$ is symmetric:} Let $\varphi_f, \varphi_g \in \mathscr{D}$, then \begin{align*}
        \langle \varphi_f, A \varphi_g\rangle &= \lim_{s\to 0} \left\langle \varphi_f, \frac{U(s)-\mathds1}{\i s}\varphi_g\right\rangle\\
        &= \lim_{s\to 0} \left\langle \frac{\mathds1-U(-s)}{\i s} \varphi_f, \varphi_g\right\rangle\\
        &= \langle -\i \varphi_{-f'}, \varphi_g \rangle \\
        &= \langle A \varphi_f, \varphi_g\rangle,
    \end{align*}
    so $A$ is symmetric.
    \item \textit{$A$ is essentially self-adjoint:} Suppose that there is a $u \in \mathscr{D}_{A^*}$ so that $A^* u = \i u$. Then for each $\varphi \in \mathscr{D}_A:=\mathscr{D}$ one gets \begin{align*}
        \frac{\dd}{\dd t} \langle U(t) \varphi, u\rangle &= \langle \i A U(t) \varphi, u\rangle\\
        &= -\i \langle U(t)\varphi, A^* u\rangle \\
        &= -\i \langle U(t) \varphi, \i u\rangle\\
        &= \langle U(t) \varphi, u\rangle.
    \end{align*}   
    Thus, $f(t)=\langle U(t) \varphi, u\rangle$ satisfies $f'=f$, so $f(t)=f(0)\e^t$. As a unitary, $U(t)$ has norm one. So $t\mapsto|f(t)|$ has to be bounded, which implies that $f(0)=\langle \varphi, u\rangle=0$. Since $\varphi \in \mathscr{D}$ was choosen arbitrarily, we can conclude that $u \equiv 0$. Similarly, one proves that $A^* u = - \i u$ has no more solutions than the trivial one. By \refcor{criterionesa}, $A$ is essentially self-adjoint on $\mathscr{D}$.
    \item \textit{Prove $U(t)=V(t)$ for $V(t)=\exp(\i t\overline{A})$:} Let $\varphi \in \mathscr{D}$. By \reflem{lemmastone} we can conclude $V(t)\varphi \in \mathscr{D}_{\overline{A}}$ and $V'(t)\varphi =\i A V(t) \varphi$ from the fact that $\varphi \in \mathscr{D}_{\overline{A}}$ which is clearly fulfilled since $\overline{A}$ is an extension of $A$. We already know that $U(t) \varphi \in \mathscr{D}\subseteq \mathscr{D}_{\overline{A}}$ for all $t\in\R$. Differentiating $w(t):=U(t)\varphi - V(t)\varphi$ yields \begin{align*}
        w'(t) &= \i A U(t) \varphi - \i \overline{A}V(t)\varphi \\
        &= \i \overline{A}w(t),
    \end{align*}
    where we replaced $A$ by $\overline{A}$ in the second step since both operators coincide on $\mathscr{D}$.\\
    Thus, by the product rule \begin{align*}
        \frac{\dd}{\dd t} \|w(t)\|^2 &=\frac{\dd}{\dd t}  \langle w(t), w(t)\rangle \\
        &= - \i \langle \overline{A} w(t), w(t)\rangle + \i \langle w(t), \overline{A}w(t)\rangle\\
        &= 0.
    \end{align*}
    Thanks to $w(0)=0$ this implies $w \equiv 0$. Hence, $U(t) \varphi = V(t) \varphi$ for all $t\in\R$, $\varphi \in \mathscr{D}$. Since $\mathscr{D}$ is dense in $\H$, we can conclude equality of $U(t)$ and $V(t)$ on $\H$.
    \item \textit{Conclude self-adjointness of $A$:} We can calculate the derivatives of $U(t)\varphi$ and $V(t)\varphi$ and evaluate them at $t=0$, respectively. Since $U(t)=V(t)$, the derivatives coincide. By \eqref{lemmastone1} this yields  \begin{align*}
       \i A\varphi = \frac{\dd}{\dd t}\Big|_{t=0} U(t) \varphi = \frac{\dd}{\dd t}\Big|_{t=0} V(t) \varphi = \i \overline{A}\varphi.
    \end{align*}
    Since this holds true for arbitrary $\varphi \in \mathscr{D}$, we can conclude self-adjointness of $A$ on $\mathscr{D}$.
    \item \textit{Uniqueness of $A$:} Suppose that there exists a self-adjoint operator $B$ such that $\exp(\i tB)=U(t)=\exp(\i tA)$. Then by \eqref{lemmastone1} we can conclude $A=B$.
\end{enumerate}

\end{proof}

\textsc{Stone's theorem} allows us to define the derivative of the mapping $t \mapsto U(t)$.

\begin{corollary}[Differential properties of semigroups]
    Let $\{U(t)\}_{t \in \R}$ be a semigroup of unitary operators with infinitesimal generator $A$. Assume $u \in \mathscr{D}_A$. Then
    \begin{enumerate}
        \item $U(t)u \in \mathscr{D}_A$ for each $t \geq 0$,
        \item $AU(t)u = U(t)Au$ for all $t > 0$,
        \item The mapping $t \mapsto U(t)$ is differentiable for each $t > 0$,
        \item $\frac{\dd}{\dd t} U(t)u = AU(t)u$.
    \end{enumerate}
    \label{differential_semigroup}
\end{corollary}

\begin{proof}
    See \cite[p.435]{Evans2010}.
\end{proof}

\subsection{Dissipativity}

The physical interpretation behind the following definition is that a dissipative operator does not increase the energy. This condition ensures that the operator either decreases or maintains the energy, which is typical in physical systems involving some dissipation effects.

\begin{definition}
    An operator $A$ on a Hilbert space $\H$ is called \textbf{dissipative} if
    \begin{align}
        \Re (\langle Ax,x \rangle ) \leq 0 \quad \forall x \in \mathscr{D}_A.
    \end{align}
    An operator $A$ is called \textbf{m-dissipative} if $A$ is dissipative and $(\mathds1-A)$ is surjective.
\end{definition}

\section{Existence results based on semigroup theory}
\subsection{The linear case}
This chapter is based on \cite[Chapter 7.5]{pazy}. Consider the linear Schrödinger equation
\begin{align*}
    \partial_t{u} = \i \Delta u -\i Vu
\end{align*}
where $V$ is a potential. We will consider the equation in the Hilbert space $\H=L^2(\R^d)$.\\

Semigroup theory provides an elegant method for constructing a solution to the initial-value problem. The right side of Schrödinger's equation will be interpreted as an operator acting on $u$. It will be shown that this operator is the infinitesimal generator of a group of unitary operators on $L^2(\R^d)$.

In the first step, define an operator $A_0$ associated with the differential operator $\i\Delta$.

\begin{definition}
    Let $\mathscr{D}_{A_0}:=H^2(\R^d)$. For $u \in \mathscr{D}_{A_0}$ we define $A_0$ by
    \begin{align*}
        A_{0}u:=\i\Delta u.
    \end{align*}
\end{definition}

Using the Fourier transform one can show that the operator $A_0$ is self-adjoint.

\begin{lemma}
    The operator $\i A:=-\Delta$ is self-adjoint in $L^2(\R^d)$.
    \label{linear_selfadjoint}
\end{lemma}

\begin{proof}
    $-\Delta$ is symmetric since
    \begin{align*}
        \langle -\Delta u,v \rangle = - \int_{\R^d} \Delta u \cdot \bar{v} \d x = - \int_{\R^d} u \cdot \overline{\Delta v} \d x = \langle u,-\Delta v \rangle
    \end{align*}
    where integration by parts was used in the second step. Let $\lambda \in \C$ with $\mathrm{Im}(\lambda) \neq 0$. To show that $-\Delta$ is self-adjoint it now suffices to show that $\overline{\mathrm{im}(\lambda I - \i A_0)}=L^2(\R^d)$. By \refcor{criterionesa} this proves essential self-adjointness of $-\Delta$ in $L^2(\R^d)$. But every self-adjoint operator admits a self-adjoint extension.
    For $f \in \mathscr{C}_\mathrm{c}^{\infty}(\R^d)$ we define
    \begin{align*}
        u(x)=(2\pi)^{-\frac{d}{2}} \int_{\R^d} \frac{\hat{f}(\eta)}{\lambda + | \eta |^2}\e^{\i x\cdot \eta} \d \eta.
    \end{align*}
    It holds $u \in H^2(\R^d)=\mathscr{D}_{A_0}$ and $u$ is the solution to $(\lambda \mathds1 - \i A_0)u=f$ since
    \begin{align*}
        (\lambda \mathds1 - \i A_0)u &= \lambda u + \Delta u = \lambda u + (2\pi)^{-\frac{d}{2}} \int_{\R^d} \frac{\hat{f}(\eta) }{\lambda + | \eta |^2}\Delta \e^{\i x\cdot\eta} \d \eta \\ &= (2\pi)^{-\frac{d}{2}} \int_{\R^d} \frac{\hat{f}(\eta) }{\lambda + | \eta |^2} (\lambda + | \eta |^2)\e^{\i x\cdot\eta} \d \eta \\ &= (2\pi)^{-\frac{d}{2}} \int_{\R^d} \hat{f}(\eta) \e^{\i x\cdot \eta} \d \eta = \check{\hat{f}}(x) = f(x) \text{ a.e.}
    \end{align*}
    Thus, $\mathscr{C}_\text{c}^{\infty}(\R^d) \subseteq \mathrm{im}(\lambda \mathds1 -\i A_0)$. Since $\mathscr{C}_\mathrm{c}^{\infty}(\R^d)$ is dense in $L^2(\R^d)$, the same holds for $\mathrm{im}(\lambda \mathds1 -\i A_0)$.
\end{proof}

Applying \textsc{Stone's theorem}, we can conclude immediately:

\begin{corollary}
    $A_0$ is the infinitesimal generator of a group of unitary operators on $L^2(\R^d)$.\label{corsemig}
\end{corollary}

In the second step, we define an operator $V$ associated with the potential $V$ by
\begin{align*}
    (Vu)(x) := V(x)u(x) \quad \text{ for all } x \in \R^d
\end{align*}
with the domain $\mathscr{D}_V=\{ u \in L^2(\R^d) : Vu \in L^2(\R^d)\}$.

We will show that $A_0-\i V$ is the infinitesimal generator of a group of unitary operators on $L^2(\R^d)$. According to \textsc{Stone's theorem}, it suffices to show that $A_0-\i V$ is self-adjoint. This will be proven using the fact that the operator is m-dissipative which follows from the perturbation theorem for m-dissipative operators.

\begin{theorem}[Perturbation theorem for m-dissipative operators]
    Let $A$ and $B$ be linear operators such that $\mathscr{D}_B \supseteq \mathscr{D}_A$ and $A+tB$ is dissipative for $0 \leq t \leq 1$. If
    \begin{align}
        \Vert Bx \Vert \leq \alpha \Vert Ax \Vert + \beta \Vert x \Vert \quad \forall x \in \mathscr{D}_A
    \end{align}
    where $0 \leq \alpha < 1, \beta \geq 0$ and $A+t_0B$ is m-dissipative for some $t_0 \in [0,1]$, then $A+tB$ is m-dissipative for all $t \in [0,1]$.
\end{theorem}

We will begin by proving the norm estimate necessary for applying the perturbation theorem. We have to restrict our considerations to the case $p > d/2$.
\begin{lemma}
    Let $V \in L^p(\R^d)$. If $p > d/2$ and $p \geq 2$, then for all $\epsilon >0$ there exists a constant $C(\epsilon) \in \R$ such that
    \begin{align}
        \Vert Vu \Vert_{L^2(\R^d)} \leq \epsilon \Vert \Delta u \Vert_{L^2(\R^d)} + C(\epsilon) \Vert u \Vert_{L^2(\R^d)}\label{pertineq}
    \end{align}
    for all $u \in H^2(\R^d)$.
\end{lemma}
\begin{proof}
    Let $u \in H^2(\R^d)$. Then $(1+|x|^2) \hat{u}(x) \in L^2(\R^d)$ by definition. Since $p > d/2$, we have $(1+|x|^2)^{-1} \in L^p(\R^d)$. By using generalized spherical coordinates we can see that
    \begin{align*}
        \int_{\R^d} (1+|x|^2)^{-p} \d x = \sigma(\mathbbm{S}^{d-1}) \int_0^\infty (1+r^2)^{-p}r^{d-1} \d r.
    \end{align*}
    As $r \rightarrow \infty$, $g(r) := \mathcal{O}(r^{-2p+d-1})$ is integrable if and only if
    \begin{align*}
        -2p+d-1 < -1 \Leftrightarrow p > \frac{d}{2}.
    \end{align*}
    By choosing appropriate conjugated exponents one can use Hölder's inequality to derive
    \begin{align*}
        \Vert fg \Vert_{L^q(\R^d)} \leq \Vert f \Vert_{L^p(\R^d)} \Vert g \Vert_{L^2(\R^d)}
    \end{align*}
    for $\frac{1}{p}+\frac{1}{2}=\frac{1}{q}$, i.e. $q=\frac{2p}{2+p}$. It follows:
    \begin{align*}
        \Vert \hat{u} \Vert_{L^q(\R^d)} &= \left( \int_{\R^d} | \hat{u}(x) |^q \d x \right)^\frac{1}{q} = \left( \int_{\R^d} \left( 1+|x|^2 \right)^{-q} \left(1+|x|^2\right)^{q} | \hat{u}(x) |^q \d x \right)^\frac{1}{q} \\ &\leq \left( \int_{\R^d} \left( 1+|x|^2 \right) ^{-p} \d x \right)^{\frac{1}{p}} \left( \int_{\R^d} \left( 1+|x|^2 \right) ^{2} | \hat{u}(x) |^2 \d x \right)^\frac{1}{2} = C_p \Vert u \Vert_{H^2(\R^d)}.
    \end{align*}
    Using the Fourier transform and \textsc{Plancherel's theorem}, one can prove the inequality
    \begin{align*}
        \Vert u \Vert_{H^2(\R^d)} \leq C\left(\Vert u \Vert_{L^2(\R^d)} + \Vert \Delta u \Vert_{L^2(\R^d)}\right)
    \end{align*}
    for some constant $C$. Inserting this yields
    \begin{align*}
         \Vert \hat{u} \Vert_{L^q(\R^d)} \leq C \left(\Vert u \Vert_{L^2(\R^d)} + \Vert \Delta u \Vert_{L^2(\R^d)}\right).
    \end{align*}
    Since $p \geq 2$, we have $1 \leq q \leq 2$ and there exists an $r \geq 2$ such that $\frac{1}{r}+\frac{1}{q}=1$. Applying the \textsc{Hausdorff-Young inequality} \eqref{hyineq} yields
    \begin{align*}
        \Vert u \Vert_{L^r(\R^d)} \leq \Vert \hat{u} \Vert_{L^q(\R^d)} \leq C \left(\Vert u \Vert_{L^2(\R^d)} + \Vert \Delta u \Vert_{L^2(\R^d)}\right).
    \end{align*}
    To achieve an arbitrarily small coefficent, replace $u(x)$ by $u(\rho x)$ for some $\rho > 0$. Then
    \begin{align*}
       \Delta u(\rho x) = \rho \sum_{i=1}^n \partial^2_i u(\rho x).
    \end{align*}
    By choosing an appropriate $\rho$, the coefficient of $\Vert \Delta u \Vert_{L^2(\R^d)}$ can be made arbitrarily small. Let $\epsilon > 0$. Choose $\rho > 0$ such that
    \begin{align*}
        \Vert u \Vert_{L^r(\R^d)} \Vert V \Vert_{L^p(\R^d)} \leq \epsilon \Vert \Delta u \Vert_{L^2(\R^d)} + C(\epsilon) \Vert u \Vert_{L^2(\R^d)}.
    \end{align*}
    Finally, since $\frac{1}{r}+\frac{1}{p}=\frac{1}{2}$ applying \textsc{Hölder's inequality} with the conjugated exponents $(\frac{1}{2}p,\frac{1}{2}r)$ yields
    \begin{align*}
        \Vert Vu \Vert_{L^2(\R^d)}^2 = \int_{\R^d} |V(x)|^2|u(x)|^2 \d x \leq \left( \int_{\R^d} |V(x)|^p \d x \right)^\frac{2}{p} \left( \int_{\R^d} |u(x)|^r \d x \right)^\frac{2}{r}
    \end{align*}
    and thus
    \begin{align*}
        \Vert Vu \Vert_{L^2(\R^d)} \leq \Vert V \Vert_{L^p(\R^d)} \Vert u \Vert_{L^r(\R^d)} \leq \epsilon \Vert \Delta u \Vert_{L^2(\R^d)} + C(\epsilon) \Vert u \Vert_{L^2(\R^d)}, 
    \end{align*}
    as desired.
\end{proof}

Now we can apply the perturbation theorem to prove that $A_0-\i V$ is self-adjoint.

\begin{theorem}
    Let $V \in L^p(\R^d,\R)$. If $p > d/2, p \geq 2$ then $A_0-\i V$ is the infinitesimal generator of a group of unitary operators on $L^2(\R^d)$.
\end{theorem}
\begin{proof}
    According to \textsc{Stone's theorem} it suffices to show that $A_0-\i V$ is self-adjoint. According to \reflem{linear_selfadjoint} $\i A_0$ is self-adjoint and $\mathrm{im}(\mathds1 \pm \i A_0)$ is dense in $L^2(\R^d)$. Thus the self-adjoint extension is surjective and $\pm A_0$ is m-dissipative. Since $V$ is a real operator, $V$ is symmetric. Thus, also $\i A_0 + V$ is symmetric. It remains to be shown that $\mathrm{im}(\mathds1 \pm (A_0 - \i V)) = L^2(\R^d)$. Since $\pm A_0$ is m-dissipative and \eqref{pertineq} holds for all $u \in \mathscr{D}_{A_0}$, the perturbation theorem states that $\pm (A_0 - \i V)$ is m-dissipative and thus, $\mathds1 \pm (A_0 - \i V)$ is surjective.
\end{proof}

Let $\{S(t)\}_{t \in \R}$ denote the semigroup of unitary operators generated by $B := A_0 - \i V$. According to \textsc{Stone's theorem} it can be written as
\begin{align*}
    S(t)=\exp(\i tB).
\end{align*}
The Schrödinger equation can be rewritten as
\begin{align*}
    \partial_t u = Bu.
\end{align*}
Define $u(t,x) := S(t)u^0(x)$. Using the differential properties of semigroups from \refcor{differential_semigroup} we can easily see that $u(t,x)$ is a solution to the initial value problem for the linear Schrödinger equation:
\begin{align*}
    \frac{\dd}{\dd t} S(t)u^0 &= BS(t)u^0,\\
    u(0,x)&=S(0)u^0(x)=\mathds1 u^0(x)=u^0(x).
\end{align*}

\subsection{Duhamel's principle}

Before analyzing the local and global existence of solutions for the nonlinear case, we will apply \textsc{Duhamel's principle} to arrive at a formula for the solution to \eqref{eq:NLS}. \textsc{Duhamel's principle} aims to write the solutions of inhomogenous linear evolution equations using the solution operator of the homogenous problem that was examined in the previous chapter.

\begin{theorem} [\textsc{Duhamel's principle}]\label{duhamel}
    Consider the linear differential operator $L:=\sum_{\alpha} c_{\alpha} D_x^{\alpha}$ on $\R^d$ where $c_{\alpha} \in \R$ are constants. Suppose that the solution to the homogenous initial value problem
    \begin{align}
         \begin{cases}
            \begin{aligned}
                \partial_t u(t,x) - Lu(t,x) &= 0,  \quad &\text{on } \R \times \R^d, \\
                u(0,x) &= u^0(x),   &\text{on } \R^d,
            \end{aligned}
        \end{cases}\label{eq:Duhamel_hom}
    \end{align}
    is given by $u(t,x)=S(t)u^0(x)$.
    Then the solution to the inhomogenous initial value problem
    \begin{align}
         \begin{cases}
            \begin{aligned}
                \partial_t u(t,x) - Lu(t,x) &= f(t,x),  \quad &\text{on } \R \times \R^d, \\
                u(0,x) &= g(x),  \quad &\text{on } \R^d ,
            \end{aligned}
        \end{cases}\label{eq:Duhamel_inhom}
    \end{align}
    is given by
    \begin{align}
        u(t,x) = S(t)g(x)+ \int_0^t S(t-s)f(s,x) \d s.
    \end{align}
\end{theorem}

We will provide a formal justification for \textsc{Duhamel's principle}. Assume $u^0, f, g$ are sufficiently smooth and we can differentiate under the integral sign without additional justification. By assumption $u_h(t,x) := S(t)u^0(x)$ satisfies the homogenous initial value problem \eqref{eq:Duhamel_hom}. Thus it holds that $u_h(x,0)=S(0)u^0(x)=\mathds1 u^0(x)=u^0(x)$ and the initial condition is satisfied:
\begin{align*}
    u(x,0) = S(0)u^0(x)+0=u^0(x).
\end{align*}
To show that $u$ satisfies the inhomogenous partial differential equation we use the differential property of the semigroup from \refcor{differential_semigroup} and compute
\begin{align*}
    \partial_t S(t)u^0(x) = LS(t)u^0(x).
\end{align*}
Using the Leibniz rule for differentiation under the integral sign yields
\begin{align*}
    \partial_t u(t,x) &= \partial_t (S(t)u^0(x)) + S(t-t)f(t,x) + \int_0^t \partial_t (S(t-s)f(s,x)) \d s \\ &= LS(t)u^0(x) + f(t,x) + \int_0^t LS(t-s)f(s,x) \d s \\ &= f(t,x) + L\left( S(t)u^0(x) + \int_0^t S(t-s)f(s,x) \d s \right) = f(t,x) + Lu(t,x).
\end{align*}
We conclude
\begin{align*}
    \partial_t u(t,x) - Lu(t,x) = f(t,x),
\end{align*}
as desired.\\

Intuitively, the inhomogenous term $f(t,x)$ can be interpreted as a source term, something external being added to the system at each time step. \textsc{Duhamel's principle} uses the superposition principle and enables us to construct a solution by superposing the effects of the forcing function. Infinitely many homogenous initival value problems are solved where $f(t,x)$ is treated as a sequence of infinitesimal impulses acting at different times $s$. Each source $f(s,x)$ is thought of as a new initial condition for a particular time step that evolves from $s$ to $t$ according to the solution of the homogenous problem. By integrating the effects over time, we arrive at the final solution.\\

Considering now the special case of the nonlinear Schrödinger equation \eqref{eq:NLS}, defining the linear operator $A_0$ on $\mathscr{D}_{A_0}$ by $A_0 u=\i \Delta u$, the initial value problem \eqref{eq:NLS} can be rewritten as 
    \begin{align*}
         \begin{cases}
            \begin{aligned}
                \partial_t u(t,x) - A_0 u(t,x) &= F(t,x),  \quad &\text{on } \R \times \R^d, \\
                u(0,x) &= u^0(x),  \quad &\text{on } \R^d ,
            \end{aligned}
        \end{cases}
    \end{align*}
where $F(t,x)= - \i \lambda |u(t,x)|^p u(t,x)$.\\
From \refcor{corsemig} we can conclude that $-A_0$ is the infinitesimal generator of a semigroup of unitary operators $\{S(t)\}_{t\in\R}$ on $L^2(\R^2)$ given by
\begin{align*}
    S(t) = \exp(A_0t)=\exp(\i t\Delta).
\end{align*}
Applying \textsc{Duhamel's principle} we can write the solution to \eqref{eq:NLS} in the following form:
\begin{align}
   \begin{split} u(t,x) &= \exp(\i t\Delta)u^0(x)+ \int_0^t \exp(\i(t-s)\Delta)F(s,x)  \d s \\ &= \exp(\i t\Delta) u^0(x) - \i \lambda \int_0^t \exp(\i(t-s)\Delta) (|u(s,x)|^p u(s,x)) \d s. \label{phiduhamel}\end{split}
\end{align}

\subsection{The nonlinear case}

This chapter is heavily based on \cite[Chapter 8.1]{pazy}, where the entire theory was developed for the case $p=2$. In the following, we generalize the results to the general case.\\
We are dealing with the nonlinear Schrödinger equation \eqref{eq:NLS} in two space dimensions, since in higher dimensions there are more subtleties regarding the use of Sobolev embeddings.
Using the Fourier transform, we can derive an explicit formula for the propagator $S(t)$. From the previous chapter we know that $u(t,x) = S(t)u^0(x)$ solves
\begin{align*}
         \begin{cases}
            \begin{aligned}
                \i \partial_t u(t,x) + \Delta u(t,x) &= 0,  \quad &\text{on } \R \times \R^d, \\
                u(0,x) &= u^0(x),  \quad &\text{on } \R^d .
            \end{aligned}
        \end{cases}
    \end{align*}

Fourier transform of $u$ in the spatial variable $x$ only transforms the equation into an ODE in $t$:
\begin{align*}
         \begin{cases}
            \begin{aligned}
                \i\partial_t \hat{u}(t,y) - | y |^2 \hat{u}(t,y) &= 0,  \quad &\text{on } \R \times \R^d, \\
                \hat{u}(0,y)&=\widehat{u^0}(y),  \quad &\text{on } \R^d .
            \end{aligned}
        \end{cases}
    \end{align*}
The solution is given by
\begin{align*}
    \hat{u}(t,y) = \e^{-\i t| y |^2} \widehat{u^0}(y).
\end{align*}
Applying the inverse Fourier transformation and using the property $\widehat{f \ast g} = (2 \pi)^{\frac{d}{2}} \hat{f} \hat{g}$ yields
\begin{align*}
    u = (2 \pi)^{-d/2} G \ast u^0,
\end{align*}
where $G$ satisfies $\hat{G} = \e^{-\i t| y |^2}$. Using the formula of the inverse Fourier transform from \refdef{Fourier} we can compute $G$ explicitly:
\begin{align*}
    G(t,x) &= \frac{1}{(2 \pi)^{\frac{d}{2}}} \int_{\R^d} \exp \left(\i  x\cdot y  - \i t | y |^2 \right) \d y = \frac{1}{(2 \i t)^{\frac{d}{2}}} \exp \left( -\frac{| x |^2}{4\i t} \right),
\end{align*}
where the complex Gaussian integral
\begin{align*}
    \int_{\R^d} \exp\left( - \frac{1}{2}  x\cdot Ax  +  b\cdot x  + c \right) \d x = \frac{(2 \pi)^\frac{d}{2}}{\sqrt{\det(A)}}\ \exp\left(\frac{1}{2}  b\cdot A^{-1}b  + c \right)
\end{align*}
for $A \in \R^{d \times d}$ symmetric positive definite was used.
For $d=2$ the solution is thus given by
\begin{align}
    (S(t)u)(x)=\frac{1}{4\pi \i t} \int_{\R^2} \exp\left(\i \frac{|x-y|^2}{4t}\right) u(y) \d y.\label{stexplicit}
\end{align}

\begin{lemma}[Decay estimate]\label{lemdecayest}
    Let $\{S(t)\}_{t\in\R}$ be a semigroup of unitary operators on $L^2(\R^2)$. If $m\in[2,\infty]$ and $\frac1m + \frac1q =1$, then $S(t)$ can be uniquely extended to a bounded operator from $L^q(\R^2)$ to $L^m(\R^2)$ and \begin{align}
        \|S(t) u\|_{L^m(\R^2)} \leq (4 \pi t)^{-(2/q-1)} \|u\|_{L^q(\R^2)}
    \end{align}
    holds.
\end{lemma}

\begin{proof}
The proof is based on a result from operator theory.

\begin{tcolorbox}\begin{center}
    Interlude: \underline{\textsc{Riesz-Thorin interpolation theorem}}
\end{center}
For a given linear operator $T$ which fulfills the inequalities \begin{align*}
    \|Tf\|_{L^{m_0}} \leq M_0 \|f\|_{L^{q_0}} \quad \text{and} \quad \|Tf\|_{L^{m_1}} \leq M_1 \|f\|_{L^{q_1}}
\end{align*}
for suitable constants $M_0, M_1 >0$, the following estimate holds: \begin{align}
    \|Tf\|_{L^m} \leq M_0^\Theta M_1^{1-\Theta} \|f\|_{L^q},
\end{align}
where $\frac1m=\frac{1-\Theta}{m_0}+\frac{\Theta}{m_1}$ and $\frac1q = \frac{1-\Theta}{q_0} + \frac{\Theta}{q_1}$.
\end{tcolorbox}
In order to be able to apply this result, we have to establish the two inequalities. 
For $(q_0,m_0)=(1,\infty)$, we use the explicit representation \eqref{stexplicit} and get \begin{align*}
    \|S(t)u\|_{L^\infty(\R^2)} \leq \underbrace{(4 \pi t)^{-1}}_{=:M_0} \|u\|_{L^1(\R^2)}.
\end{align*}
To obtain the second inequality, we choose $(q_1,m_1)=(2,2)$. In this case, $S(t)$ is an isometry, that is, \begin{align*}
    \|S(t)u\|_{L^2(\R^2)}=\|u\|_{L^2(\R^2)}
\end{align*}
holds for $u \in L^2(\R^2)$. So the second inequality is actually an equality with $M_1:=1$. \\
We end up with \begin{align*}
    \|S(t)u\|_{L^m(\R^2)} \leq [(4\pi t)^{-1}]^\Theta \|u\|_{L^q(\R^2)},
\end{align*}
where $\Theta = 2/q -1$.

\end{proof}

The strategy is as usual: we use Lipschitz estimates to prove local existence. The proof is based on \refthm{uniquethm}, and in \reflem{lemmaestimatescont}, the necessary Lipschitz estimates are established.

\begin{theorem}[Uniqueness theorem]\label{uniquethm}
    Let $Y$ be a Banach space and let $F:[t_0,T] \times Y \to Y$ be uniformly Lipschitz continuous in $Y$ and for each $y \in Y$ let $F(\cdot,y):[t_0,T] \to Y$ be continuous. If $u^0 \in \mathscr{D}_A$, where $-A$ is the generator of a semigroup $\{S(t)\}_{t\in\R}$, then the initial value problem 
    \begin{align}
        \begin{cases}
            \frac{\dd u}{\dd t}(t)+Au(t)=F(t,u(t)), & t>t_0,\\
            u(t_0) = u^0, &
        \end{cases}\label{ivpsemi}
    \end{align}
    has a unique classical solution on $[t_0,T]$.
    \end{theorem}

\begin{proof}
    See \cite[Theorem 6.1.7.]{pazy}.
\end{proof}

\begin{lemma}\label{lemmaestimatescont}
    Let $p \geq 1$. The nonlinear mapping $F(u):= \i \lambda |u|^p u$ maps $H^2(\R^2)$ into itself and satisfies for $u,v \in H^2(\R^2)$:
    \begin{enumerate}
        \item $\|F(u)\|_{H^2(\R^2)} \leq C \|u\|_{L^\infty(\R^2)}^p \|u\|_{H^2(\R^2)}$ for a $C>0$,\label{lemmanonlinear1}
        \item $\|F(u)-F(v)\|_{H^2(\R^2)} \leq C\left( \|u\|_{H^2(\R^2)}^p + \|v\|_{H^2(\R^2)}^p\right) \|u-v\|_{H^2(\R^2)}$ for a $C >0$.\label{lemmanonlinear2}
    \end{enumerate}
\end{lemma}

\begin{proof}
    From the Sobolev embedding theorem in $\R^2$, we have $H^2(\R^2) \hookrightarrow L^\infty(\R^2)$ (cf. \eqref{sobemb1}). That is, for $u\in H^2(\R^2)$ there is a $C>0$ s.t. \begin{align*}
        \|u\|_{L^\infty(\R^2)} \leq C \|u\|_{H^2(\R^2)}.
    \end{align*}
    This ensures that $F$ is well-defined and bounded in $L^\infty(\R^2)$ since the following holds for every $u \in L^\infty(\R^2)$ \begin{align*}
        \||u|^p u\|_{L^\infty(\R^2)}= \sup_{x \in \R^2} |u(x)|^{p+1} \leq M^{p+1},
    \end{align*}
    where $M>0$.\\
    To prove \ref{lemmanonlinear1}, we calculate the derivatives of $|u|^p u$ up to second order:
\begin{align}
    \D(|u|^p u)&= p |u|^{p-2} u (\D u)u+ |u|^p \D u = (p+1) |u|^p \D u, \\
    \D^2(|u|^p u) &= p |u|^{p-2} u \left(u \D^2 u + (p-1) |\D u|^2\right).\label{seconderivativelem}
\end{align}

    The $H^2$-norm of $F$ is given by
    \begin{align*}
        \|F(u)\|_{H^2(\R^2)} = \sqrt{\|F(u)\|_{L^2(\R^2)}^2 + \|\D F(u)\|_{L^2(\R^2)}^2 + \|\D^2 F(u)\|_{L^2(\R^2)}^2}.
    \end{align*}
    Now we want to find estimates on the three $L^2$-norms:
    \begin{itemize}
        \item \underline{Bound on $\|F(u)\|_{L^2(\R^2)}$:}
        \begin{align*}
    \|F(u)\|_{L^2(\R^2)} = |\lambda| \||u|^p u\|_{L^2(\R^2)} \leq |\lambda| \| u\|_{L^\infty(\R^2)}^p \|u\|_{L^2(\R^2)}.
    \end{align*}
        \item \underline{Bound on $\|\D F(u)\|_{L^2(\R^2)}$:} 
        \begin{align*}
            \|\D F(u)\|_{L^2(\R^2)} = \|\i \lambda (p+1) |u|^p \D u\|_{L^2(\R^2)}
            &\leq |\lambda|(p+1) \|u\|^p_{L^\infty(\R^2)} \|\D u\|_{L^2(\R^2)} \\
            &\leq |\lambda|(p+1) \|u\|^p_{L^\infty(\R^2)} \|u\|_{H^2(\R^2)}.
        \end{align*}
        \item \underline{Bound on $\|\D^2 F(u)\|_{L^2(\R^2)}$:} 
        We have \begin{align*}
            \||u|^p \D^2 u\|_{L^2(\R^2)} \leq \|u\|_{L^\infty(\R^2)}^p \|\D^2 u\|_{L^2(\R^2)} \leq \|u\|_{L^\infty(\R^2)}^p \|u\|_{H^2(\R^2)},
        \end{align*} 
        which provides a bound on the first term of the second derivative \eqref{seconderivativelem}. For the second term we use a \textsc{Gagliardo-Nirenberg inequality} of the form \begin{align*}
            \|\D u\|_{L^4(\R^2)} \leq C \|u\|^{1/2}_{L^\infty(\R^2)}\|u\|^{1/2}_{H^2(\R^2)}.
        \end{align*}
        This leads to \begin{align*}
            \| |u|^{p-1} |\D u|^2\|_{L^2(\R^2)} \leq \|u\|^{p-1}_{L^\infty(\R^2)} \||\D u|^2\|_{L^2(\R^2)}&=\|u\|^{p-1}_{L^\infty(\R^2)} \|\D u\|_{L^\infty(\R^4)}  \\
            &\leq \|u\|^{p-1}_{L^\infty(\R^2)} \left( C \|u\|_{L^\infty(\R^2)}^{1/2} \|u\|_{H^2(\R^2)}^{1/2}\right)^2\\
            &= C^2\|u\|_{L^\infty(\R^2)}^{p}\|u\|_{H^2(\R^2)}.
        \end{align*}
        We end up with \begin{align*}
            \|\D^2 F(u)\|_{L^2(\R^2)} \leq |\lambda| (p+ (p-1)C^2)\|u\|_{L^\infty(\R^2)}^{p}\|u\|_{H^2(\R^2)}.
        \end{align*}
    \end{itemize}
The combination of the three bounds yields \begin{align*}
    \|F(u)\|_{H^2(\R^2)} \leq |\lambda|\max\{1, p+1, p+(p-1)C^2\} \|u\|_{L^\infty(\R^2)}^{p}\|u\|_{H^2(\R^2)} =: C \|u\|_{L^\infty(\R^2)}^{p}\|u\|_{H^2(\R^2)}.
\end{align*}

    The Lipschitz estimate \eqref{lemmanonlinear2} is a straight calculation. We decompose the difference $F(u)-F(v)=\i \lambda (|u|^p u - |v|^p v)$ as \begin{align*}
        |u|^p u - |v|^p v= G(u,v)(u-v),
    \end{align*}
where $G(u,v):=\sum_{j=0}^p |u|^{p-j} |v|^j$. Clearly, one has $|G(u,v)|\leq (p+1)\max\{|u|^p,|v|^p\}$. This yields \begin{align*}
    \|F(u)-F(v)\|_{H^2(\R^2)} = |\lambda| \||u|^p u - |v|^p v\|_{H^2(\R^2)} \leq |\lambda| (p+1)\max\{\|u\|_{L^\infty(\R^2)}^p,\|v\|_{L^\infty(\R^2)}^p\} \|u-v\|_{H^2(\R^2)}.
\end{align*}
The usual Sobolev embedding inequality $\|u\|_{L^\infty(\R^2)} \leq \Tilde C \|u\|_{H^2(\R^2)} $ leads to the desired result with $C:=|\lambda| \Tilde{C}^p$:
\begin{align*}
    \|F(u)-F(v)\|_{H^2(\R^2)} \leq C \left( \|u\|_{H^2(\R^2)}^p + \|v\|_{H^2(\R^2)}^p\right) \|u-v\|_{H^2(\R^2)}.
\end{align*}
\end{proof}

\begin{corollary}[Local existence of a solution]\label{corlocexistence}
    For every $u^0 \in H^2(\R^2)$ there exists a unique solution $u$ to \eqref{ivpsemi} defined for $t \in [0,T_\mathrm{max})$ such that \begin{align*}
        u \in \mathscr{C}^1([0,T_\mathrm{max}),L^2(\R^2)) \cap \mathscr{C}([0,T_\mathrm{max}),H^2(\R^2))
    \end{align*}
    with the property that either $T_\mathrm{max}=\infty$ or $T_\mathrm{max}<\infty$ and $\lim_{t \to T_\mathrm{max}} \|u\|_{H^2(\R^2)}=\infty$.
\end{corollary}
\begin{proof}
From \refthm{uniquethm} we obtain local existence of a unique classical solution since $F(u):=|u|^p u$ is locally Lipschitz as a map $H^2(\R^2) \to H^2(\R^2)$ by \reflem{lemmaestimatescont}. Hence, the initial value problem \eqref{eq:NLS} has a unique solution $u \in \mathscr{C}([0,T_\text{max}),H^2(\R^2))$. To gain higher regularity, we make use of a classical bootstrap argument: since the equation has the form $\partial_t u = -A u + F(u(t))$ where $A=-\Delta$ with domain $\mathscr{D}_A=H^2(\R^2)$ and since $u(t)\in H^2(\R^2)$ for all $t \in [0,T_\text{max})$, it follows that both $Au(t)\in L^2(\R^2)$ and $F(u(t))\in H^2(\R^2) \subseteq L^2(\R^2)$. Thus the right-hand side of the equation lies in $L^2(\R^2)$ und is continuous. Consequently, $u \in \mathscr{C}^1([0,T_\text{max}),L^2(\R^2))$.\\
    In order to prove uniqueness also in the weak sense, we proceed as follows: Let $u_1, u_2$ be two solutions with the same initial condition. We define their difference \begin{align*}
        w(t):=u_1(t)-u_2(t).
    \end{align*}
    Subtracting the equations for $u_1$ and $u_2$, we obtain \begin{align}
        \i \partial_t w + \Delta w = |u_1|^p u_1 - |u_2|^p u_2.\label{hilfsgleichungw}
    \end{align}
    Taking the $H^2$-norm, we compute \begin{align*}
        \frac{\dd}{\dd t} \|w(t)\|_{H^2(\R^2)}^2 = 2 \Re \langle \partial_t w, w\rangle.
    \end{align*}
    Using \eqref{hilfsgleichungw} we obtain \begin{align*}
        \frac{\dd}{\dd t}\|w(t)\|_{H^2(\R^2)}^2 &= 2 \Re \langle -\i (\Delta w + |u_1|^p u_1 - |u_2|^p u_2), w\rangle\\
        &= 2 \Re \langle -\i ( |u_1|^p u_1 - |u_2|^p u_2), w\rangle,
    \end{align*}
    where the last equality follows from the fact that $\langle -\i \Delta w, w\rangle$ is purely imaginary.\\
    The application of \reflem{lemmaestimatescont} to the difference $F(u_1)-F(u_2)$ yields \begin{align*}
        \frac{\dd}{\dd t}\|w(t)\|_{H^2(\R^2)}^2 \leq C \left(\|u_1\|_{H^2(\R^2)}^p + \|u_2\|_{H^2(\R^2)}^p\right) \|w(t)\|_{H^2(\R^2)}^2.
    \end{align*}
    Now we can use \textsc{Gronwall's inequality}. If $\frac{\dd}{\dd t} f(t) \leq g(t) f(t)$, then \begin{align*}
        f(t) \leq f(0) \exp\left(\int_0^t g(s) \d s \right).
    \end{align*}
    Setting $f(t)=\|w(t)\|_{H^2(\R^2)}^2$ and $g(t)=C \left(\|u_1\|_{H^2(\R^2)}^p + \|u_2\|_{H^2(\R^2)}^p\right) $, we obtain \begin{align*}
    \|w(t)\|_{H^2(\R^2)}^2 \leq \|w(0)\|_{H^2(\R^2)}^2 \exp\left( C\int_0^t \left(\|u_1\|_{H^2(\R^2)}^p + \|u_2\|_{H^2(\R^2)}^p\right) \d s \right).
    \end{align*}
    Thanks to $u\in \mathscr{C}([0,T_\text{max}),H^2(\R^2))$, the integral is finite. Since $w(0)=0$, it follows that $\|w(t)\|_{H^2(\R^2)}^2=0$ for all $t$, which implies $w(t) \equiv 0$.\\
    In order to discuss the blow-up conditions, we investigate the equation \begin{align*}
        \frac{\dd}{\dd t}\|u(t)\|_{H^2(\R^2)}^2 = 2 \Re \langle \partial_t u, u \rangle = 2 \Re \langle -\i(\Delta u - |u|^p u), u \rangle.
    \end{align*}

Since the contributions of the Laplacian term as well as of the nonlinear term vanishes (both inner products are purely imaginary), the derivative of the $H^2$-norm of $u(t)$ is zero, so $\|u(t)\|_{H^2(\R^2)}$ is constant, so in particular bounded, which guarantees local existence. Hence, the solution cannot blow up unless $\|u^0\|_{H^2(\R^2)}$ was already infinite, which contradicts the assumption that $u^0 \in H^2(\R^2)$. Thus, if $T_\text{max}$ is finite while $\|u(t)\|_{H^2(\R^2)}$ remained bounded, we could extend the solution past $T_\text{max}$, contradicting its maximality. The only way $T_\text{max}$ can be finite is if the $H^2$-norm of the solution diverges: 
\begin{align*}
    \lim_{t \to T_\text{max}} \|u\|_{H^2(\R^2)}=\infty.
\end{align*}
Therefore, the only way a solution can cease to exist in finite time is through a norm blow-up.

\end{proof}

\begin{lemma}\label{lemglobsol}
    Let $u^0 \in H^2(\R^2)$ and let $u$ be a solution to \eqref{eq:NLS} on $[0,T)$. If $\lambda \geq 0$ then $\|u(t)\|_{H^2(\R^2)}$ is bounded on $[0,T)$.
\end{lemma}

\begin{proof}
    The proof consists of three steps.
    \begin{enumerate}
        \item \textit{$L^2(\R^2)$-norm preservation:}\label{lemstep1} We multiply \eqref{eq:NLS} by $\overline{u}$ and integrate over $\R^2$:
        \begin{align}
             \int_{\R^2} \left( \frac1{\i}\frac{\partial u}{\partial t} \overline{u} - (\Delta u) \overline{u} + \lambda |u|^{p+2}\right) \d x=0.\label{lemmaintnonl}
        \end{align}
        Using $\frac{\dd}{\dd t} |u|^2 = \frac{\partial u}{\partial t} \overline{u} + u \frac{\partial \overline{u}}{\partial t}$ we get \begin{align*}
            \frac{\dd}{\dd t} \int_{\R^2} |u|^2 \d x = 2 \Re \int_{\R^2} \frac{\partial u}{\partial t} \overline{u}  \d x.
        \end{align*}
        Therefore, the first term in \eqref{lemmaintnonl} is given by \begin{align}
            \int_{\R^2}  \frac1{\i}\frac{\partial u}{\partial t} \overline{u} \d x = \frac{1}{2\i} \frac{\dd}{\dd t} \int_{\R^2} |u|^2 \d x.\label{lemnonint1}
        \end{align}
        The second term in \eqref{lemmaintnonl} can be treated by using integration by parts:
        \begin{align}
            \int_{\R^2} - \Delta u \cdot \overline{u} \d x = \int_{\R^2} \nabla u \cdot \nabla \overline{u} \d x = \int_{\R^2} |\nabla u|^2 \d x.\label{lemnonint2}
        \end{align}
    Now we take the imaginary part of all of the three terms. Since the expression in \eqref{lemnonint2} and also the integral of the third term in \eqref{lemmaintnonl} are purely real, the imaginary part equals zero. Therefore, only imaginary part of \eqref{lemnonint1} is relevant and we obtain \begin{align*}
        \frac{1}{2} \frac{\dd}{\dd t} \int_{\R^2} |u|^2 \d x=0.
    \end{align*}
    That is, \begin{align*}
\|u(t)\|_{L^2(\R^2)}=\|u^0\|_{L^2(\R^2)}.
    \end{align*}
    \item\textit{Boundedness of $\|u(t)\|_{H^1(\R^2)}$:}\label{lemstep2}
    We multiply \eqref{eq:NLS}  by $\frac{\partial \overline{u}}{\partial t}$:
    \begin{align*}
             \int_{\R^2} \left( \frac1{\i}\frac{\partial u}{\partial t}  - \Delta u + \lambda |u|^{p}u\right)\frac{\partial \overline{u}}{\partial t} \d x=0.
    \end{align*}
    We want to take the real part of this integral to conclude the asserted statement. The first term vanishes:
    \begin{align*}
        \Re \left( \int_{\R^2}  \frac1{\i}\frac{\partial u}{\partial t}  \frac{\partial \overline{u}}{\partial t} \d x\right) = \Re \left( \int_{\R^2}  \frac1{\i}\left|\frac{\partial u}{\partial t}\right|^2   \d x\right) =0.
    \end{align*}
    Using $\frac{\dd}{\dd t}|\nabla u|^2=2\Re\left(\nabla u \cdot \nabla \frac{\partial \overline{u}}{\partial t}\right)$ we can compute the real part of the second term, using integration by parts again:
    \begin{align*}
        \Re\left(\int_{\R^2}(-\Delta u) \frac{\partial \overline{u}}{\partial t} \d x \right) = \Re\left(\int_{\R^2}\nabla u \cdot \nabla \frac{\partial \overline{u}}{\partial t} \d x \right) = \frac12 \frac{\dd}{\dd t}\int_{\R^2} |\nabla u|^2 \d x.
    \end{align*}
    The real part of the third, nonlinear term can be computed as \begin{align*}
        \Re\left(\int_{\R^2}\lambda|u|^p u \frac{\partial \overline{u}}{\partial t} \d x \right) = \frac{\dd}{\dd t} \int_{\R^2} \frac{\lambda}{p+2} |u|^{p+2} \d x.
    \end{align*}
    This can be done by calculating $\frac{\dd}{\dd t} \int_{\R^2} G(u) \d x $ for $G(u)=\frac{\lambda}{p+2} |u|^{p+2}$.\\
    Putting all terms together yields \begin{align*}
        \frac{\dd}{\dd t}  \int_{\R^2} |\nabla u|^2 \d x + \frac{\dd}{\dd t} \int_{\R^2} \frac{\lambda}{p+2} |u|^{p+2} \d x =0.
    \end{align*}
    Hence, the energy \begin{align*}
        E(t)=\int_{\R^2}  |\nabla u|^2 \d x + \int_{\R^2} \frac{\lambda}{p+2} |u|^{p+2} \d x
    \end{align*} is conserved, i.e. $E(t)=E(0)$ for all $t \in [0,T)$. Thanks to non-negativity of $E$, the $L^2$-norm of $\nabla u(t)$ is bounded for all $t \in [0,T)$. From step \ref{lemstep1} we already know that $\|u(t)\|_{L^2(\R^2)}$ is also bounded, which leads to \begin{align*}
        \|u(t)\|_{H^1(\R^2)} = \sqrt{\|u(t)\|_{L^2(\R^2)}^2 + \|\nabla u(t)\|_{L^2(\R^2)}^2}\leq C
    \end{align*} for $C>0$, as claimed.
    \item\textit{Boundedness of $\|u(t)\|_{H^2(\R^2)}$:} Thanks to \eqref{phiduhamel}, the solution $u$ of the integral equation ($S(\cdot)$ denotes the time evolution from the previous chapter on the linear Schrödinger equation and $F$ is the inhomogenity $\lambda |u|^p u$) is given by \begin{align*}
        u(t)= S(t) u^0 - \int_0^t S(t-s) F(u(s)) \d s.
    \end{align*}
    We apply any first order derivative $\D$ to obtain \begin{align*}
        \D u(t) = S(t) \D u^0 - \int_0^t S(t-s) \D F(u(s)) \d s.
    \end{align*}
    The derivative of the nonlinearity is given by $\D F(u)=\lambda u^p \D u + \lambda p |u|^{p-2} u \Re(\overline{u} \D u)$. It is sufficient to control the contribution of the first term since the second term is of the same order due to $|\Re(\overline{u}\D u)| \leq |u| |\D u|$.\\
    We fix the Sobolev exponent $m>2$ and choose the corresponding conjugate exponents:
    \begin{align*}
        q= \frac{m}{m-1} \quad \text{and} \quad r=\frac{2m}{m-2}.
    \end{align*}
    Applying Hölder's inequality yields \begin{align*}
        \| |u|^p \D u\|_{L^q(\R^2)} \leq \|u\|_{L^r(\R^2)}^p \|\D u \|_{L^2(\R^2)}.
    \end{align*}
    Using the Sobolev embedding $H^1(\R^2) \hookrightarrow L^r(\R^2)$ (cf. \eqref{sobemb2}) we have \begin{align*}
    \|u\|_{L^r(\R^2)} \leq C \|u\|_{H^1(\R^2)},
    \end{align*}
    which is bounded thanks to step \ref{lemstep2}. Since $\|\D u \|_{L^2(\R^2)} \leq C \|u\|_{H^1(R^2)}$, it follows that \begin{align*}
        \| |u|^p \D u \|_{L^q(\R^2)} \leq C.
    \end{align*}

Using the semigroup estimate from \reflem{lemdecayest}, the semigroup $S(t-s)$ maps $L^q$ to $L^m$ with decay:
\begin{align*}
    \|S(t-s) \D F(u(s))\|_{L^m(\R^2)} \leq C (t-s)^{-(2/q-1)} \|\D F(u(s))\|_{L^q(\R^2)}.
\end{align*}
Substitute $q=\frac{m}{m-1}$ and compute the exponent \begin{align*}
    2/q -1 = \frac{2(m-1)}{m}-1 = \frac{m-2}{2}.
\end{align*}
Thus, the time decay becomes $(t-s)^{-(m-2)/m}=(t-s)^{-1+2/m}$. The integral \begin{align*}
    \int_0^t (t-s)^{-1+2/m} \d s
\end{align*}
converges because $-1+2/m > -1$ for $m>2$. This gives: \begin{align*}
    \|\D u(t)\|_{L^m(\R^2)} \leq C \|u^0\|_{H^2(\R^2)} + C(T).
\end{align*}

By Sobolev embedding $W^{1,m}(\R^2) \hookrightarrow L^\infty(\R^2)$ for $m>2$ we have \begin{align*}
    \|u(t)\|_{L^\infty(\R^2)} \leq C.
\end{align*}
Finally, we use \textsc{Gronwall's inequality} on the $H^2(\R^2)$-norm:
\begin{align*}
    \|u(t)\|_{H^2(\R^2)} \leq \|u^0\|_{H^2(\R^2)}+ C \int_0^t \|u(s)\|^p_{L^\infty(\R^2)} \|u(s)\|_{H^2(\R^2)} \d s.
\end{align*}
Since $\|u(s)\|_{L^\infty(\R^2)} \leq C$, Gronwall implies: \begin{align*}
    \|u(t)\|_{H^2(\R^2)} \leq C \|u^0\|_{H^2(\R^2)} \exp(C T).
\end{align*}
Thus, $\|u(t)\|_{H^2(\R^2)}$ is bounded on $[0,T)$.
     \end{enumerate}
\end{proof}

\begin{theorem}[Global existence of a solution]\label{globexsg}
Let $u^0 \in H^2(\R^2)$. If $\lambda \geq 0$ then the initial value problem \eqref{eq:NLS} has a unique global solution 
\begin{align*}
        u \in \mathscr{C}^1([0,\infty),L^2(\R^2)) \cap \mathscr{C}([0,\infty),H^2(\R^2)).
    \end{align*}
\end{theorem}

\begin{proof}
    By \refcor{corlocexistence}, the existence of a unique local solution on a time interval $[0,T_\text{max})$ is guaranteed. Moreover, we can conclude from $T_\text{max}<\infty$ that $\lim_{t \to T_\text{max}} \|u(t)\|_{H^2(\R^2)}=\infty$.\\ \reflem{lemglobsol} ensures that if $\lambda \geq 0$, the $H^2(\R^2)$-norm of $u$ is bounded on $[0,T)$. That is, there is a $T$-dependent constant $C(T)$ such that $\|u(t)\|_{H^2(\R^2)} \leq C(T)$ for all $t \in [0,T)$. Via a continuation argument this holds for any $T>0$, so we can conclude $\lim_{t \to T_\text{max}} \|u(t)\|_{H^2(\R^2)} < \infty$.\\ Combining these facts, we must have $T_\text{max}=\infty$. So the local solution exists for all $t\geq 0$ and is unique by \refcor{corlocexistence}. It follows that \eqref{eq:NLS} is globally well-posed, as claimed. \end{proof}
\section{Existence results based on Strichartz estimates}
Strichartz estimates provide valuable insights into the behaviour of solutions to NSE. Some inspirations for the proof structures presented in the following are taken from \cite{courant}. 
\subsection{Strichartz spaces}

\begin{definition}[Strichartz space]
Let $\Omega_t \subseteq \R$ and $\Omega_x\subseteq \R^d$. For each measurable function $f:\Omega_t \times \Omega_x \to \C$ we define
      \begin{align*}
      \|f\|_{L_t^{q} L_x^{r}(\Omega_t \times \Omega_x)} = \left( \int_{\Omega_t} \left( \int_{\Omega_x} |f|^{r} \d x \right)^{\frac{q}{r}} \d t \right)^{\frac{1}{q}}.
  \end{align*}
  The normed space \begin{align*}
      L_t^{q} L_x^{r}(\Omega_t \times \Omega_x):= \left\{ f:\Omega_t \times \Omega_x \to \C \text{ measurable } \Big| \|f\|_{L_t^{q} L_x^{r}(\Omega_t \times \Omega_x)} < \infty \right\}
  \end{align*}
  is called \textbf{Strichartz space} or \textbf{mixed Lebesgue space}.
\end{definition}

The following theorem is standard and can be found in \cite[Theorem 2.3]{Tao2006}.

\begin{lemma}
Fix $d\geq 1$ and call a pair $(q,r)$ of exponents \textbf{admissible} if $2 \leq q,r \leq \infty$, $\frac2q + \frac{d}{r}=\frac{d}{2}$ and $(q,r,d)\neq(2,\infty,2)$. Then for any admissible exponents $(q,r)$ and $(\Tilde{q},\Tilde{r})$ we have the homogeneous Strichartz estimate \begin{align}
    \| \exp(\i t \Delta) u^0\|_{L_t^q L_x^r(\R \times \R^d)} \leq C_{d,q,r} \|u^0\|_{L_x^2(\R^d)}.\label{Strich1}
\end{align}
Further suppose that $(\Tilde{q}, \Tilde{r})$ satisfy the same restrictions as $(q,r)$ and $(\Tilde{q}', \Tilde{r}')$ are their Hölder conjugates, then the dual homogeneous Strichartz estimates take the form
\begin{align}
    \left\| \int_\R \exp(-\i s \Delta) F(s) \d s \right\|_{L_x^2(\R^d)} \leq C_{d, \Tilde{q}, \Tilde{r}} \|F\|_{L_t^{\Tilde{q}'}L_x^{\Tilde{r}'}(\R \times \R^d)}.\label{Strich2}
\end{align}
The inhomogeneous (or retarded) Strichartz estimates are \begin{align}
    \left\|\int_{-\infty}^t \exp(\i(t-s) \Delta) F(s) \d s \right\|_{L_t^q L_x^r(\R \times \R^d)} \leq C_{d,q,r,\Tilde{q}, \Tilde{r}} \|F\|_{L_t^{\Tilde{q}'}L_x^{\Tilde{r}'}(\R \times \R^d)}.\label{Strich3}
\end{align}
\end{lemma}

\subsection{Local existence of solutions}

From now on, we consider the space
\begin{align*}
X := \mathscr{C}([0,T],L^2(\R^d)) \cap L_t^q L_x^r([0,T] \times \R^d)
\end{align*}
for suitable Strichartz pairs $(q, r)$. \\
By choosing $T$ sufficiently small, we can ensure that the solution to the integral equation \eqref{phiduhamel} provided by \textsc{Duhamel's principle} is a contraction mapping. By \textsc{Banach's fixed point theorem}, this mapping has a unique fixed point in $X$, giving a local solution to \eqref{eq:NLS}.

\begin{theorem}
    Let $X := \mathscr{C}([0,T],L^2(\R^d)) \cap L_t^q L_x^r([0,T] \times \R^d)$, where $(q,r)$ is a Strichartz-admissible pair. Choose $T$ maximal such that
    \begin{align}
        T < (4 C |\lambda|R^p)^{-\frac{1}{1-(p+q)/2}} \quad \text{ and } \quad C T^{p/q} \left( \|u\|_{L_t^\infty L_x^r}^p + \|v\|_{L_t^\infty L_x^r}^p \right) < \frac{1}{2}
    \end{align}
    for all $u, v \in X$ where $C$ is the constant from the Strichartz estimate.
    Then the map 
    \begin{align}
        \Phi : X \to X, \quad u \mapsto \Phi(u) := \exp(\i t \Delta) u^0 - \i \lambda \int_0^t \exp(\i (t-s) \Delta) \big( |u(s)|^p u(s) \big) \, \mathrm{d}s
    \end{align}
    defines a contraction on a ball in $X$ with radius $R :=4C \|u^0\|_{L_x^2(\R^d)}$.
    \label{localexistencestrichartz}
\end{theorem}

Before proving this statement, we provide two useful lemmata.

\begin{lemma}
Let $u:[0,T]\times \R^d \to \C$ be a local solution to the IVP \eqref{eq:NLS}, provided it exists. Let $R :=4C \|u^0\|_{L_x^2(\R^d)}$. Then the following inequality holds true:\label{lem:hell}
\begin{align}
    \||u|^p u\|_{L_t^{\tilde{q}'} L_x^{\tilde{r}'}}\leq R^{p+1} T^{1-\frac{p+1}{q}}.\label{boundhell}
\end{align}
    $\tilde{q}'$ and $q$ are related by $\frac{1}{\tilde{q}'}+\frac1q=1$.
\end{lemma}

\begin{proof}
    First, we prove an ancillary estimate $
        \||u|^p u\|_{L_t^{\tilde{q}'}L_x^{\tilde{r}'}} \leq \|u\|_{L_t^q L_x^{r(p+1)}}^{p+1}$ by using \textsc{Hölder's inequality}. Expanding the $ L_t^{\tilde{q}'} L_x^{\tilde{r}'} $-norm of $ |u|^p u$ yields
  \begin{align*}
      \||u|^p u\|_{L_t^{\tilde{q}'} L_x^{\tilde{r}'}([0,T] \times \R^d)} = \left( \int_0^T \left( \int_{\R^d} |u(t, x)|^{(p+1) \tilde{r}'} \, \d x \right)^{\frac{\tilde{q}'}{\tilde{r}'}} \d t \right)^{\frac{1}{\tilde{q}'}}.
  \end{align*}
   To express this norm in terms of $\| u \|_{L_t^q L_x^{r(p+1)}}$, we use \textsc{Hölder's inequality} with the exponents chosen so that $\tilde{q}'$ and $\tilde{r}'$ are conjugate to $q$ and $r(p+1)$, respectively:
   \begin{align*}
      \frac{1}{\tilde{q}'} + \frac{1}{q} = 1 \quad \text{and} \quad \frac{1}{\tilde{r}'} + \frac{1}{r(p+1)} = 1. 
   \end{align*}

   We apply \textsc{Hölder’s inequality} in the $x$-variable first:
   \begin{align*}
       \int_{\R^d} |u(t, x)|^{(p+1) \tilde{r}'} \d x \leq \left( \int_{\R^d} |u(t, x)|^{r(p+1)} \d x \right)^{\frac{\tilde{r}'}{r(p+1)}}.
   \end{align*}
Raising both sides to the power $\frac{\tilde{q}'}{\tilde{r}'}$, we obtain
 \begin{align*}
    \left( \int_{\R^d} |u(t, x)|^{(p+1) \tilde{r}'} \d x \right)^{\frac{\tilde{q}'}{\tilde{r}'}} \leq \left( \int_{\R^d} |u(t, x)|^{r(p+1)} \d x \right)^{\frac{\tilde{q}'}{r(p+1)}}.    
 \end{align*}

   Next, we integrate in $t$ and apply \textsc{Hölder’s inequality} again, now in the $t$-variable:
   \begin{align*}
         \int_0^T \left( \int_{\R^d} |u(t, x)|^{r(p+1)} \d x \right)^{\frac{\tilde{q}'}{r(p+1)}} \d t \leq \left( \int_0^T \left( \int_{\R^d} |u(t, x)|^{r(p+1)} \d x \right)^{\frac{q}{r(p+1)}} \d t \right)^{\frac{\tilde{q}'}{q}}.
   \end{align*}

   Combining the inequalities from the spatial and temporal applications of \textsc{Hölder’s inequality}, we have:
   \begin{align*}
      \||u|^p u\|_{L_t^{\tilde{q}'} L_x^{\tilde{r}'}([0,T] \times \R^d)} \leq \| u \|_{L_t^q L_x^{r(p+1)}([0,T] \times \R^d)}^{p+1},
   \end{align*}
   as desired.\\
Using \textsc{Hölder's inequality} with $\frac{1}{a} + \frac{1}{b} = 1$, we have:
    \begin{align*}
        \int_0^T f(t) \, \mathrm{d}t \leq \left( \int_0^T |f(t)|^a \, \mathrm{d}t \right)^\frac{1}{a} \left( \int_0^T 1^b \, \mathrm{d}t \right)^\frac{1}{b} = T^\frac{1}{b} \left( \int_0^T |f(t)|^a \, \mathrm{d}t \right)^\frac{1}{a}.
    \end{align*}
To estimate $\|u\|_{L_t^q L_x^{r(p+1)}}^{p+1}$, set $a = \frac{q}{p+1}$ and $b = \frac{q}{q - (p+1)}$:
    \begin{align*}
        \|u\|_{L_t^q L_x^{r(p+1)}}^{p+1} \leq T^\frac{q-(p+1)}{q} \left( \int_0^T \|u(t, \cdot)\|_{L_x^{r(p+1)}}^q \, \mathrm{d}t \right)^\frac{p+1}{q}.
    \end{align*}
    Since $r(p+1) > r$, one has $\|u(t, \cdot)\|_{L_x^{r(p+1)}} \leq \|u(t, \cdot)\|_{L_x^r}$(\footnote{This holds true since we consider only the local problem.}). Then, for $u \in \B_R(u^0)$, we have the bound
    \begin{align*}
        \int_0^T \|u(t, \cdot)\|_{L_x^{r(p+1)}}^q \, \mathrm{d}t \leq \int_0^T \|u(t, \cdot)\|_{L_x^r}^q \, \mathrm{d}t = \|u\|_{L_t^q L_x^r}^q \leq R^q.
    \end{align*}
    Thus,
    \begin{align*}
        \|u\|_{L_t^q L_x^{r(p+1)}}^{p+1} \leq R^{p+1} T^{1 - \frac{p+1}{q}}.
    \end{align*}
\end{proof}

\begin{lemma} Let $u,v:[0,T]\times \R^d \to \C$ be two local solution to the IVP \eqref{eq:NLS}, provided they exist. Then the following estimate holds:\label{lemestimatephi} 
    \begin{align}
        \|\Phi(u) - \Phi(v)\|_X \leq C \left( \|u\|_{L_t^q L_x^r}^p + \|v\|_{L_t^q L_x^r}^p \right) \|u - v\|_X.\label{estimatephi}
    \end{align}
\end{lemma}

\begin{proof}
    We begin by writing the difference: \begin{align*}
        \Phi(u)-\Phi(v)= -\i \lambda \int_0^t \exp(\i(t-s)\Delta) \left(|u(s)|^p u(s)-|v(s)|^p v(s)\right) \d s.
    \end{align*}
    Applying the Strichartz estimate \eqref{Strich3}, we obtain \begin{align*}
        \|\Phi(u)-\Phi(v)\|_X \leq C \| |u|^p u -|v|^p v \|_{L_t^{q'} L_x^{r'}}\leq C \|(|u|^p + |v|^p) |u - v|\|_{L_t^{q'} L_x^{r'}},
    \end{align*}
    where the elementary inequality $\big| |u|^p u - |v|^p v \big| \leq C (|u|^p + |v|^p) |u - v|$ was used in the last step.\\
    No we apply \textsc{Hölder's inequality} in time and space. Since $u,v \in L_t^{q} L_x^{r}$, we have $f:=|u|^p + |v|^p \in L_t^{q/p} L_x^{r/p}$ and $g:= u-v \in L_t^{q} L_x^{r}$. The Hölder exponents are determined by 
    \begin{align*}
        \frac1{q'}=\frac{1}{q/p}+\frac1q \quad \text{and} \quad \frac1{r'}=\frac{1}{r/p}+\frac1r.
    \end{align*}
    \textsc{Hölder's inequality} $\|fg\|_{L_t^{q'} L_x^{r'}} \leq \|f\|_{L_t^{q/p} L_x^{r/p}} \|g\|_{L_t^{q} L_x^{r}}$ yields \begin{align*}
        \|(|u|^p + |v|^p) |u - v|\|_{L_t^{q'} L_x^{r'}} \leq \||u|^p + |v|^p\|_{L_t^{q/p} L_x^{r/p}} \|u - v\|_{L_t^{q} L_x^{r}} \leq \left(\|u\|_{L_t^{q} L_x^{r}}^p + \|v\|_{L_t^{q} L_x^{r}}^p\right)\|u - v\|_{L_t^{q} L_x^{r}}.
    \end{align*}
    We obtain by definition of the norm on $X$:
    \begin{align*}
        \|\Phi(u)-\Phi(v)\|_X \leq  \left(\|u\|_{L_t^{q} L_x^{r}}^p + \|v\|_{L_t^{q} L_x^{r}}^p\right)\|u - v\|_{L_t^{q} L_x^{r}} \leq \left(\|u\|_{L_t^{q} L_x^{r}}^p + \|v\|_{L_t^{q} L_x^{r}}^p\right)\|u - v\|_X.
    \end{align*}
\end{proof}

\begin{proof}[Proof of \refthm{localexistencestrichartz}]
    First, we show that $\Phi$ maps a small ball around the initial condition into itself. Define the ball $\B_R(u^0) := \{ u \in X : \|u - u^0 \|_X \leq R \}$, where the $X$-norm is defined by $\|u\|_X := \|u\|_{L_t^\infty L_x^2} +\|u\|_{L_t^q L_x^r}$.

    Using the homogeneous Strichartz estimate \eqref{Strich1}, we get
    \begin{align*}
        \|\exp(\i t \Delta) u^0\|_{L_t^q L_x^r} \leq C \|u^0\|_{L_x^2}.
    \end{align*}
    This ensures that the contribution of the free evolution stays close to the initial condition $u^0$ in the $L_t^q L_x^r$-norm. The inhomogeneous Strichartz estimate \eqref{Strich3} allows us to bound the nonlinear part of $\Phi$:
    \begin{align*}
        \left\| \int_0^t \exp(\i (t-s) \Delta) \big( |u(s)|^p u(s) \big) \, \mathrm{d}s \right\|_{L_t^q L_x^r} \leq C_{d,q,r,\tilde{q}, \tilde{r}} |\lambda| \, \||u|^p u\|_{L_t^{\tilde{q}'} L_x^{\tilde{r}'}}.
    \end{align*}

Combining the estimates for the linear and nonlinear parts gives
\begin{align*}
        \|\Phi(u)\|_X =& \|\Phi(u)\|_{L_t^\infty L_x^2} +\|\Phi(u)\|_{L_t^q L_x^r}\\
         \leq& \|\exp(\i t \Delta) u^0\|_{L_t^\infty L_x^2} + \left\| \int_0^t \exp(\i (t-s) \Delta) \big( |u(s)|^p u(s) \big) \, \mathrm{d}s \right\|_{L_t^\infty L_x^2} \\
         &+\|\exp(\i t \Delta) u^0\|_{L_t^q L_x^r} + \left\| \int_0^t \exp(\i (t-s) \Delta) \big( |u(s)|^p u(s) \big) \, \mathrm{d}s \right\|_{L_t^q L_x^r} \\
         \leq& C \|u^0\|_{L_x^2} + C |\lambda| \||u|^p u\|_{L_t^{\tilde{q}'} L_x^{\tilde{r}'}} + C \|u^0\|_{L_x^2} + C |\lambda| \||u|^p u\|_{L_t^{\tilde{q}'} L_x^{\tilde{r}'}}\\
         =& 2C \|u^0\|_{L_x^2} + 2 C |\lambda|\||u|^p u\|_{L_t^{\tilde{q}'} L_x^{\tilde{r}'}}.
    \end{align*}
Now we can use our assumption on $R$ and get $2C\|u^0\|_{L_x^2}=\frac{R}{2}$. The second term can be bounded by the bound \eqref{boundhell} from \reflem{lem:hell}:
\begin{align*}
    2 C |\lambda|\||u|^p u\|_{L_t^{\tilde{q}'} L_x^{\tilde{r}'}} \leq 2 C |\lambda| R^{p+1} T^{1-\frac{p+1}{q}}.
\end{align*}
Due to the definition of $T$, this is bounded by $R/2$, too. This yields \begin{align*}
    \|\Phi(u)\|_X \leq R,
\end{align*}
as claimed.

    It remains to show that $\Phi$ is Lipschitz continuous. Here we use the bound on $\|\Phi(u) - \Phi(v)\|_X$ from \reflem{lemestimatephi}.

    Applying \textsc{Hölder's inequality} in time, we have $\|f\|_{L_t^q} \leq T^{1/q} \|f\|_{L_t^\infty}$. Thus,
    \begin{align*}
        \|u\|_{L_t^q L_x^r}^p + \|v\|_{L_t^q L_x^r}^p \leq T^{p/q} \left( \|u\|_{L_t^\infty L_x^r}^p + \|v\|_{L_t^\infty L_x^r}^p \right).
    \end{align*}
    Since $T$ is, by definition, small enough so that $C T^{p/q} \left( \|u\|_{L_t^\infty L_x^r}^p + \|v\|_{L_t^\infty L_x^r}^p \right) < \frac{1}{2}$, we achieve Lipschitz continuity with modulus $\frac{1}{2}$. The result then follows by \textsc{Banach's fixed-point theorem}.
\end{proof}

\section{Discussion}

In this work, we explored two different mathematical approaches to studying the local existence of solutions to the nonlinear Schrödinger equation, namely by semigroup theory arguments and by using Strichartz estimates. Both methods provide valuable insights, yet they differ in elegance, computational complexity, and their connection to functional analysis. \\

Semigroup theory offers an elegant framework deeply embedded in functional analysis and operator theory. It is based on the interpretation of the given PDE as an operator equation -- usually an integral equation. A key feature of this method is that it accomodates a broad class of evolution equations beyond Schrödinger-type problems; for instance the heat equation and other evolution equations. We can always interpret the evolution of the system as the action of an operator semigroup, making the method more general and adaptable to various settings. While the class of equations semigroup theory can be applied to is broad, the use of Sobolev embedding theorems pose a challenge in higher dimensions.  We thus restricted ourselves to the dimension $d=2$. Local existence and uniqueness of solutions relies on \refthm{uniquethm}. Application of this uniqueness theorem required local Lipschitz continuity of the nonlinear mapping $F(u):= \i\lambda |u|^pu$. The local solution is either defined for $[0, \infty)$ and thus global or a norm blow up happens after finite time $T_\mathrm{max} < \infty$. Global existence of solutions is shown for the special case of the defocusing NSE for $\lambda > 0$ already discussed in the introduction. Here, the norm of a local solution stays bounded in $H^2$ on $[0,T)$ and a norm blowup is prevented. Thus, the solution exists globally. Global existence for the focusing NSE with $\lambda < 0$ can be proven for $p=2$ under the additional condition that
\begin{align*}
    | \lambda | \Vert u^0 \Vert_{L^2(\R^2)} < 2,
\end{align*}
hence the $L^2$-norm of the initial datum has to be sufficiently small (cf. \cite[p.233]{pazy}). This condition ensures that $\Vert u \Vert_{H^2(\R^2)}$ stays bounded. 

On the other hand, Strichartz estimates were specifically developed for Schrödinger operators and provide an alternative strategy for establishing the well-posedness of solutions. The proof strategy using Strichartz estimates might be  more intuitive and accessible as it is reminiscent of the proof of \textsc{Picard-Lindelöf's Theorem}. There the differential equation is first transformed into an integral equation, then \textsc{Banach's Fixed-point theorem} is applied to prove the existence of a solution. The reformulation of the NSE was achieved by using \textsc{Duhamel's principle} which allowed the definition of a corresponding mapping $\Phi(u)$. Similar to the proof of \textsc{Picard-Lindelöf's Theorem} it was shown that $\Phi$ is a contraction if the domain of $u$ is chosen small enough. The approach is based on the definition of an appropriate function space where the Lipschitz continuity of the nonlinear operator could be proven, akin to the strategy used in the semigroup ansatz. 

Strichartz estimates enabled us to show local existence for arbitrary dimensions $d$ which the semigroup theory only provided for $d=2$ due to the fact that Sobolev embeddings are more complicated in higher dimensions.

Overall, our analysis highlights the strengths and limitations of both methods. The semigroup approach provides a conceptually elegant and generalizable framework applicable to a broad class of evolution equations, while the Strichartz estimates approach is specifically tailored to dispersive equations like the Schrödinger equation. While the generalization to higher dimensions by semigroup theory requires significantly more effort, Strichartz estimates provide local existence in arbitrary dimensions. Future work could explore hybrid techniques that combine the advantages of both approaches to yield even more refined existence results for nonlinear Schrödinger equations and related PDEs.


\vspace{0.3cm}

\noindent {\it Acknowledgments.} 
We would like to thank Daniel Matthes for his professional and personal support in carrying out and completing this project.

\newpage\appendix
\color{black}

\newpage
\bibliographystyle{plain}  
\bibliography{bibliography}   

\begin{thebibliography}{10}

\bibitem{water_waves}
Georgi Gary Rozenman \& Lev Shemer \&~Ady Arie.
\newblock Observation of accelerating solitary wavepackets.
\newblock {\em Physical Review E}, 2020.

\bibitem{Brezis}
Haim Brezis.
\newblock {\em Functional Analysis, Sobolev Spaces and Partial Differential
  Equations}.
\newblock Universitext. Springer, 2011.

\bibitem{courant}
Thierry Cazenave.
\newblock {\em Semilinear Schrödinger Equations}.
\newblock American Mathematical Society, 2003.

\bibitem{Evans2010}
Lawrence~C Evans.
\newblock {\em Partial Differential Equations}.
\newblock Graduate studies in mathematics. American Mathematical Society,
  Providence, RI, 2 edition, March 2010.

\bibitem{BEC}
Hichem Hajaiej \& Slim Ibrahim \&~Nader Masmoudi.
\newblock Ground state solutions of the complex gross pitaevskii equation
  associated to exciton-polariton bose-einstein condensates.
\newblock {\em Journal de Mathématiques Pures et Appliquées}, 148, 2021.

\bibitem{pazy}
Amnon Pazy.
\newblock {\em Semigroups of Linear Operators and Applications to Partial
  Differential Equations}.
\newblock Applied Mathematical Sciences. Springer, Providence, RI, 1 edition,
  1983.

\bibitem{Simon}
Michael Reed \&~Barry Simon.
\newblock {\em Functional Analysis}.
\newblock Modern Methods of Mathematical Physics. American Mathematical
  Society, 1980.

\bibitem{Strichartzwave}
Robert Strichartz.
\newblock Restriction of fourier transform to quadratic surfaces and decay of
  solutions of wave equations.
\newblock {\em Duke Math. J.}, 1977.

\bibitem{Tao2006}
Terence Tao.
\newblock {\em Nonlinear Dispersive Equations: Local and Global Analysis},
  volume 106 of {\em CBMS Regional Conference Series in Mathematics}.
\newblock American Mathematical Society, 2006.

\bibitem{Blowup}
Terence Tao \& Monica Visan \&~Xiaoyi Zhang.
\newblock Minimal-mass blowup solutions of the mass-critical nls.
\newblock {\em arXiv preprint}, math/0609690, 2006.

\end{thebibliography}

\end{document}